\documentclass[24pt]{article}
\usepackage{amssymb}  
\usepackage{amsfonts}
\usepackage{amsmath}
\usepackage[T1]{fontenc}
\usepackage[italian]{babel}
\usepackage{pifont}
\usepackage{caption}
\usepackage{placeins}
\usepackage{pstricks}
\usepackage{pdfpages}
\usepackage{amsfonts}
\usepackage{graphicx}
\usepackage{placeins}
\usepackage{amsmath}
\usepackage{authblk}
\usepackage{pstricks}
\usepackage{pstricks-add}
\usepackage{framed}
\usepackage{bm}
\usepackage[utf8]{inputenc}

\usepackage{aurical}
\usepackage[T1]{fontenc}

\usepackage{xcolor,colortbl}

\usepackage{calligra}

\usepackage{tikz,fp,ifthen,fullpage}
\usepackage{pgfmath}
\usetikzlibrary{backgrounds}
\usetikzlibrary{decorations.pathmorphing,backgrounds,fit,calc,through}
\usetikzlibrary{arrows}
\usetikzlibrary{shapes,decorations,shadows}
\usetikzlibrary{fadings}
\usetikzlibrary{patterns}
\usetikzlibrary{mindmap}
\usetikzlibrary{decorations.text}
\usetikzlibrary{decorations.shapes}

\newtheorem{propr}{Propriet\`a}

\begin{document}

\bibliographystyle{plain}

\title{La matematica armonia dei suoni naturali. Ovvero, l'armonica matematica dei suoni naturali.}

\author{\it F.~Talamucci}
\affil{\footnotesize\it DIMAI, Dipartimento di Matematica e Informatica ``Ulisse Dini'',\\
\footnotesize\it Universit\`a degli Studi di Firenze, Italy\\
\it e-mail: federico.talamucci@unifi.it}
\date{}

\maketitle

\begin{abstract}
	
	\noindent
	La collaborazione della matematica nei due sistemi musicali dei suoni pitagorici e dei suoni equabili \`e evidente ed opportuna per generare gli elementi e gestirne le relazioni.
	L'unica nozione indispensabile per un intervento razionale ai due mondi sonori \`e quella di distanza fra suoni.
	La scala dei suoni naturali non sembra dialogare in modo diretto con un modulo puramente matematico che articoli la sua definizione ed \`e generalmente presentata come aggiustamento manuale degli antichi suoni pitagorici.
	Il modello di proporzioni che ha formato ed ispirato l'architettura del Cinquecento, epoca in cui la scala naturale emerge, fornisce un lucido spunto per 
	assestare in modo sistematico i suoni naturali; dall'altra parte, la generazione di suoni basata ed avviata esclusivamente sul rispetto delle proporzioni propone sistemi di suoni naturali interessanti.
	L'oggetto non \`e tanto quello di illustrare il congegno armonico della gi\`a formata scala naturale, quanto attuare la formazione di elementi sonori che rispondano al principio unitario della proporzione e che si combinino armonicamente.

\end{abstract}

\section{Preambolo}
La musica \`e nata da un impulso emotivo; laddove l'ispirazione e l'istinto sono sorretti dal pensiero 
e da fondamenti teorici si assiste ad un'affermazione e ad un progresso talmente formidabili da avere pochi eguali nelle conquiste dell'attivit\`a umana e nell'equilibrio tra natura ed intelletto. Questo \`e avvenuto per la nostra cosiddetta musica della cultura occidentale, a partire dal mondo classico fino ai nostri giorni.

\noindent
Nella decifrazione e l'elaborazione del pensiero musicale interviene certamente anche la scienza, talvolta la matematica: quanto alla presenza di quest'ultima, seppur ogni tanto forzata o azzardata o ingombrante, ritengo di poter ricondurre idealmente il punto di partenza alla geniale sintesi del pensiero musicale antico della Figura 1 e rintracciare nel corso della storia due snodi impareggiabili: i rapporti ben proporzionati fra suoni che ispirano i moduli architettonici del Rinascimento e la presa di coscienza del 
suono come sovrapposizione di suoni pi\`u semplici in epoca illuministica.

\begin{figure}
	\includegraphics[width=0.36\textwidth]{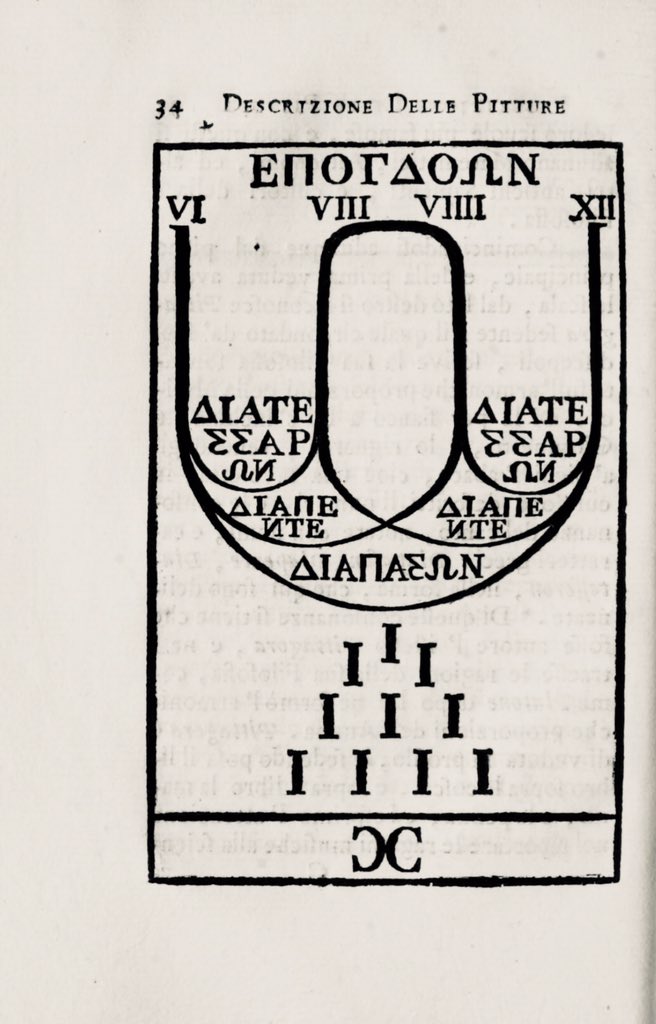}
	\caption{Prospetto della suddivione del diapason, secondo la teoria musicale pitagorica. L'illustrazione \`e tratta da \cite{bellori}.}
\end{figure}

\noindent
Riguardo al secondo episodio, che non tratteremo direttamente, va almeno ricordata nella prima met\`a del 700 la sorprendente concordanza fra la teoria scientifica dei suoni armonici\footnote{La prima dimostrazione sperimentale circa l'esistenza dei suoni armonici, che concorrono con frequenze differenti a formare un unico suono, \`e ad opera di J.~Sauveur e risale al 1701.} nella corda vibrante e la coeva teoria musicale, presente quest'ultima principalmente attraverso i trattati di J.~P.~Rameau, il primo dei quali (\cite{rameau}) gi\`a attraverso le prime parole ``{\it La musica \`e la scienza dei suoni...''} assegna all'armonia il ruolo di scienza deduttiva, su cui compiere dimostrazioni e sintesi matematiche. 
L'ispirazione \`e la medesima dei {\it Principia} di Newton nel campo della scienza del moto, di non tanto precedenti; l'impatto culturale che non permetter\`a agli studiosi dei rispettivi generi di evitare da l\`\i$\,$ in poi di averne a che fare, il medesimo.

\noindent
Il clima evocato trova una rapida sintesi in questo: pur nella percezione di un suono netto, la natura assegna al mezzo in vibrazione un insieme di suoni che si fondono; la base naturale dell'armonia musicale, che dispone i suoni e forma gli accordi, non pu\`o non rivolgersi a tale evidenza 
per formare le regole.

\noindent
L'altra circostanza di mediazione matematica in materia di musica dapprima nominata trae anch'essa ispirazione dall'insegnamento della natura: l'architetto del Rinascimento \`e prima di tutto uno scienziato che ha il compito di apprendere dall'Universo l'armonico sistema di rapporti matematici per edificare mediante esso.
Ed \`e qui che interviene la musica: i rapporti universalmente validi vengono rivelati dai suoni, che risultano in armonia se disposti nelle giuste proporzioni. Ovvero, nella musica la natura espone distintamente i rapporti valevoli in assoluto.

\noindent
La convinzione di un ordine universale della Natura filtrato dalla musica attraversa tutta l'Architettura del 400 e del 500, da Leon Battista Alberti\footnote{Condensando lo sterminato argomento in sbrigativa aneddotica citiamo la celebre e significativa ammonizione dell'Alberti, durante la costruzione del Tempio Malatestiano, di non alterare le proporzioni dei pilastri altrimenti ``{\it si discorda tutta quella musica''.}}(1404--1472) ad Andrea Palladio\footnote{Col medesimo animo della precedente nota, rievochiamo una memoria al Palladio richiesta nel 1567 sulle nuove costruzioni di Brescia in cui si legge “...le proporzioni delle voci sono armonia delle orecchie, cos\`\i$\,$ quelle delle misure sono armonia degli occhi nostri...”.} (1508--1580). L'armonia dell'Universo si rivela in musica attraverso un rapporto ben proporzionato tra i suoni: accogliere i rapporti musicali in architettura mediante il linguaggio universale dei numeri conferisce il giusto compimento al progetto da realizzare.

\noindent
Pi\`u che un connubio, una simbiosi fra Architettura rinascimentale e musica \`e un contemplarsi ed ispirarsi su piani separati, il pi\`u elevato dei quali \`e abitato dalla musica con il suo ruolo speciale e sovrano di custodire e dispensare le proporzioni universali.
Ben inteso, nessuna architettura \`e possibile senza proporzioni metriche, anche nel Medioevo o in epoca moderna se ne fa uso, ma il rivolgersi alla teoria musicale per procedere \`e un tratto distintivo del Rinascimento. Talvolta si ha l'impressione di un affidamento, un'aspettativa verso la musica da parte di alcuni umanisti pi\`u autentica e risoluta dei musicisti stessi, come pure di un'incrollabile convinzione delle norme assolute della musica da parte degli architetti rinascimentali che non si pu\`o 
non trarre ispirazione dalle elaborazioni di costoro.

\noindent
Una sorta di ``collocazione spaziale'' dei suoni mediante i rapporti \`e comunque antica almeno quanto l'Armonia delle sfere della metafisica pitagorica: il modello dei corpi celesti su sfere ruotanti \`e
un concerto nello spazio pi\`u che nel tempo, pi\`u statico che dinamico, pi\`u volto a contemplare l'eufonica concordanza dei rapporti che il piacevole fluire della musica.

\noindent
\`E proprio il mondo antico di Pitagora che ci riaccosta all'emblematica Figura 1\footnote{Una raffigurazione di contenuto pressoch\'e identico appare nella tavoletta innanzi a Pitagora, nel  celebre affresco ``{\it La Scuola di Atene}'' eseguito da Raffaello intorno al 1510. La riproduzione qui utilizzata, tratta da \cite{bellori},  appartiene a G.~P.~Bellori (1613--1696), celebre biografo di pittori seicenteschi e profondo studioso dell'arte di Raffaello e della filosofia di Platone.}, nella quale sono disposti in modo essenziale i soggetti musicali e le proporzioni aritmetiche a cui dedicare l'attenzione e le precisazioni della prossima Sezione, senza rinunciare alla praticit\`a di affiancare all'antica simbologia la trasposizione del punto di vista moderno.

\section{Suoni ed intervalli}

\noindent
La percezione acuto--grave di un suono, ovvero l'{\it altezza}, \`e correlata al numero di vibrazioni che compie un mezzo per produrre un suono. Riferendoci all'esperienza della corda vibrante, esercitata sin dall'antichit\`a, la maggiore lunghezza produce un suono pi\`u grave, la diminuzione di questa d\`a origine ad un suono pi\`u acuto.
Ai nostri giorni \`e semplice formulare in termini assoluti questa correlazione adoperando la {\it frequenza}, ovvero il numero di vibrazioni al secondo: \`e infatti ben nota la {\it proporzionalit\`a inversa tra la lunghezza della corda e la frequenza di vibrazione}, propriet\`a che viene dedotta dall'olandese Isaac Beeckman con una dimostrazione matematica solo intorno al 1614, ed approfondita negli aspetti matematici e sperimentali da Marin Mersenne (\cite{mersenne}), che prende in esame anche la tensione e la massa della corda.

\noindent
Prima di allora, la non identificata frequenza non ha comunque impedito di fissare una scala dell'altezza del suono, per lo meno relativa, basandosi su aspetti pi\`u facilmente ispezionabili della corda, come appunto la lunghezza. Ci\`o che ha colmato l'assenza di esperienza fisica del suono 
\`e una valutazione di carattere percettivo, dunque non scientifica ma inattaccabile: raddoppiando o dimezzando la lunghezza della corda si avverte una sorta di {\it replica}, pi\`u acuta o pi\`u grave, del suono prodotto dalla corda di partenza, nel senso di risultare fortemente assimilabile, estremamente consonante.

\noindent
La certezza di riascoltare pressoch\'e il medesimo suono dividendo in due la corda fa decifrare parte della riga in alto con i numeri romani ${\rm VI}$ (sei) e ${\rm XII}$ (dodici) della tavoletta, come pure la scritta $\Delta {\rm IA}\Pi {\rm A}\Sigma \Omega {\rm N}$ ({\it diapason}, ovvero ``{\it attraverso tutti''}) che abbraccia i numeri: attribuendo al suono di base il numero sei\footnote{La scelta, consolidata dalla presenza in testi importanti come i Commenti che Marsilio Ficino aggiunse alla traduzione latina del {\it Timeo} (\cite{timeo}), il dialogo platonico delle armonie e dei medi proporzionali, agevola la produzione di esempi numerici di proporzioni, come pure rileviamo nei Trattati di Architettura del Cinquecento.}, si oltrepassano tutti i suoni per ritrovarlo in corrispondenza della corda dimezzata, ma contrassegnato dal 
{\it numero doppio} dodici per porre maggior rilevanza sulla dislocazione pi\`u acuta.
Non siamo affatto lontani dal preciso concetto di frequenza, tanto che la proporzionalit\`a inversa di quest'ultima alla lunghezza della corda \`e azzardatamente talvolta attribuita a Pitagora.
D'altronde, ``vedere'' i suoni come relazione tra numeri \`e pienamente nello spirito pitagorico: la realt\`a \`e conoscibile perch\'e misurabile.

\noindent
Avvalendosi per comodit\`a della nozione di frequenza, accreditiamo il giusto valore al tema che \`e emerso: l'unica certezza universalmente riconosciuta da qualunque cultura musicale ed in ogni epoca storica, evitando di andare alla ricerca di incerte eccezioni, \`e proprio il {\it percepire e riconoscere la massima consonanza raddoppiando o dimezzando il numero di vibrazioni}\footnote{Nel gi\`a nominato Trattato \cite{rameau}, di taglio newtoniano, l'asserzione assume il ruolo di un Principio ({\it ``Principe de l'identit\'e des octaves}'') da cui si deduce matematicamente l'armonia}.

\noindent
Attraversando la gamma dei suoni udibili si scoprono dunque degli elementi ricorsivi che semplificano decisamente la gestione musicale del complesso sonoro: i suoni sono raggruppabili in una sorta di ripartizione in classi ciascuna delle quali contenente un suono con frequenza $\nu$ e quelli con frequenza $2^k \nu$, $k$ numero intero\footnote{La denominazione stessa dei suoni a cui siamo abituati assegna il medesimo nome alla classe, non al singolo suono: ad esempio con 440 vibrazioni al secondo il suono si chiama LA, con 880 vibrazioni \`e il LA pi\`u acuto, a 110 vibrazioni il LA pi\`u grave. 
Nell'ambito dell'udibilit\`a, all'incirca da 20 a 20.000 vibrazioni al secondo, i LA udibili sono complessivamente dieci, dalla frequenza 27,5 fino a 14.080, i primi otto dei quali fanno parte della tastiera di un pianoforte moderno.
Va detto che l'esatta identificazione nota--frequenza era ignota nel passato e proviene dai progressi della moderna metrologia: \`e singolare il fatto che il LA a 440, suono di riferimento per l'intonazione degli strumenti musicali, viene ufficialmente fissato dall'Organizzazione Internazionale per la Normazione I.~S.~O.~solo nell'Ottobre del 1939, seppur prendendo atto di una prassi generalizzata e diffusa ormai da tempo.}.

\noindent
Il ``ripetersi'' del suono alle frequenze 
$\nu$, $2\nu$, $\dots, \; 2^n \nu$ qualunque sia la frequenza $\nu$ -- 
nonch\'e il passo gi\`a compiuto di assegnare la comune dicitura {\it diapason} a tutti gli intervalli tra un suono ed il suo doppio --
propone, anzi impone una disposizione equidistanziata di tali suoni: 
questo avviene se la distanza fra due suoni  -- o meglio, in termini musicali, dell'{\it intervallo tra due suoni}) viene quantificata tramite il {\it rapporto delle frequenze} --  ponendo al numeratore la frequenza del suono pi\`u acuto;
in questo modo, l'intervallo {\it diapason} misura sempre $2$, da abbinare all'$1/2$ del dimezzamento della corda. La misura dell'intervallo come il rapporto tra le frequenze (oppure con quantit\`a identificative analoghe, come nel caso della tavoletta) predispone il sostegno adatto per operare con i suoni: se un suono ha altezza $\mu$ compresa nel diapason fra $\nu$ e $2\nu$, la replica $2\mu$ ha medesima distanza dagli estremi del diapason di competenza $2\nu$, $4\nu$, dato che $2\mu:2\nu=\mu:\nu$. Ovvero, in ciascun diapason si ripetono inalterate le medesime distanze tra suoni omologhi.
Altres\`\i, la misura come rapporto impone una regola di tipo moltiplicativo per la composizione di intervalli: se infatti un terzo suono ha frequenza $\eta>\mu$, la separazione  
$\dfrac{\eta}{\nu}= \dfrac{\eta}{\mu}\times \dfrac{\mu}{\nu} \times$ fa s\`\i$\;$ che
la ripartizione dell'intervallo da $\nu$ a $\eta$ nei due intervalli contigui da $\nu$ a $\mu$ e da $\mu$ a $\eta$ abbia come misura il prodotto delle misure, qualunque siano le altezze $\nu$, $\mu$ e $\eta$. 

\noindent
Una seconda ricognizione alla tavoletta della Figura 1 permette l'accesso ad ulteriori decifrazioni, ampliando con coerenza gli elementi forniti: lo spazio denominato ${\rm \Delta I A}\Pi {\rm ENTE}$ (diapente), traducibile come ``{\it attraverso cinque}'', collega i numeri ${\rm VI}$ (sei) e ${\rm VIIII}$ (nove), come pure ${\rm VIII}$ (otto) e ${\rm XII}$ (dodici), sottintendendo una porzione di $6/9=8/12=2/3$ della corda iniziale; d'altra parte, la scritta ${\rm \Delta IATE\Sigma \Sigma AP\Omega N}$ (diatessaron), ``{\it attraverso quattro}'', congiunge ${\rm VI}$ e ${\rm VIII}$, oltre che ${\rm VIIII}$ e ${\rm XII}$, presupponendo una frazione di corda pari a $6/8=9/12=3/4$ di quella iniziale. L'antico diagramma non riporta solo il nome dello spazio che separa due suoni, ma prosegue ed esaurisce la dichiarazione delle consonanze: il suono prodotto dalla corda di lunghezza $2/3$ e quello prodotto dalla corda di lunghezza $3/4$ (rispetto alla lunghezza della corda iniziale) sono in buon accordo, in armonia con il suono di partenza. L'accrescimento musicale consiste nell'acquisire la consonanza dell'intervallo {\it diapente} corrispondente alla {\it quinta} nell'armonia musicale moderna, di lunghezza (coerentemente con l'enunciata regola del rapporto) $9/6=12/8=3/2$ e dell'intervallo {\it diatessaron},  l'odierna {\it quarta}, di lunghezza $8/6=12/9=4/3$.

\noindent
La parte in basso della tavoletta, con i primi quattro numeri naturali disposti a piramide, 
riproduce la celebre {\it tetraktys} di Pitagora\footnote{Nella dottrina pitagorica, in cui il fondamento numerico stabilisce ogni aspetto della realt\`a, il complesso 1, 2, 3 e 4 conserva un significato simbolico e filosofico superiore ed \`e alla base dell'essenza dell'ordine dell'Universo.}, 
la cui presenza da una parte sta ad indicare la pratica generazione dei tre intervalli descritti con la divisione della corda in due, tre, quattro parti\footnote{Precisamente, ciascuna delle due parti uguali di suddivione produce il diapason, la divisione in tre parti genera il diapente ed il diapente pi\`u acuto, nelle rispettive lunghezze di $2/3$ e $1/3$; la ripartizione in quattro parti d\`a origine al diapason del diapason ($1/4$ della lunghezza), diapason ($1/2$), diatessaron ($3/4$).}, dall'altra appare come fermo monito a fermarsi qui nella ricerca di consonanze. \`E come se i quattro supremi numeri che si affacciano in basso risuonassero la loro autorit\`a nel perfetto ed esclusivo intreccio di loro elaborazioni che appare in alto. La {\it tetraktys} prescrive le armonie, ma \`e la musica a dare una conferma, una riprova alla filosofia: nel globale concetto pitagorico del Numero come fondamento dell'Universo la sacralit\`a dei rapporti numerici 1:2, 2:3 e 3:4 prende voce e trova riscontro nella meravigliosa sintonia dei corrispondenti suoni. 

\subsection{I suoni pitagorici} 

\noindent
Manca per ora la spiegazione dell'iscrizione in alto e della doppia ${\rm C}$ speculare in basso: quest'ultima \`e in realt\`a una ${\rm X}$ (come appare pi\`u chiaramente nell'affresco di Raffaello) a testimonianza della sacra decade 10, somma dei numeri che compongono la {\it tetraktys} e simbolo dell'Universo.

\noindent
Pi\`u importante \`e la scritta della parte superiore ${\rm E\Pi O \Gamma \Delta O \Omega N}$ ({\it epogdoon}) che incombe per grandezza ma appare esclusa dalle consonanze: la traduzione ``{\it sopra l'ottavo}'' ricapitola l'operazione di sommare un ottavo al numero otto, per ottenere $9/8$, coerentemente con gli estremi ${\rm VIII}$ e ${\rm VIIII}$ che abbraccia la scritta, a distanza appunto $9/8$. Al tempo stesso, l'epogdoon proviene dagli intervalli consonanti, come intervallo da aggiungere al diatessaron per ottenere il diapente, in accordo con la regola moltiplicativa $\frac{4}{3}\times \frac{9}{8}=\frac{3}{2}$. 

\noindent
L'ubicazione privilegiata e per niente sfuggevole dell'epogdoon sembra suggerire la rilevanza e la 
gravosit\`a della proposta di un progetto, quello di aumentare i suoni, quello di formare scale musicali.
La prospettiva \`e nuova: non si \`e aggiunta una consonanza, ma un intervallo pi\`u piccolo dei precedenti, con cui infittire il diapason, aggiungendo nuovi suoni ai soli quattro finora presenti e rendendo il diapente e il diatessaron ottemperanti a ci\`o che dichiarano, ovvero: 
quali ``cinque'' oppure ``quattro'' suoni attraversano?

\noindent
La direzione prevedibile per andare alla ricerca di altri suoni, senza trasgredire le consonanze dichiarate, \`e quella di formare l'intervallo diapente con i suoni presenti: l'unico che d\`a origine ad una novit\`a (a meno di diapason, tramite raddoppi o dimezzamenti) \`e il suono {\rm VIIII}, il cui diapente va a fermarsi sul suono di altezza 27 che oltrepassa il diapason e che corrisponde alla replica (pi\`u bassa di due diapason) 27/4, stavolta all'interno del diapason VI--XII. Si prosegue in modo analogo\footnote{Questo modo di procedere non \`e unico, ad esempio c'\`e a disposizione anche la consonanza del diatessaron; tuttavia, i vari risultati ottenuti essenzialmente si equivalgono.}
generando ogni volta la consonanza diapente dal nuovo suono per trovare (compiendo i necessari dimezzamenti) $\frac{27}{4}\rightarrow \frac{81}{8}\rightarrow\frac{243}{32}\rightarrow \frac{729}{64}$.
La collocazione dei suoni ottenuti nel diapason ``normalizzato'' tra 1 e 2 espone i suoni in modo pi\`u agevole ed \`e conveniente per le comparazioni successive: la divisione per 6 produce la sequenza -- in ordine crescente ed includendo naturalmente i gi\`a esistenti suoni $\frac{1}{6}\cdot 8 = \frac{3}{2}$ e $\frac{1}{6}\cdot 8 = \frac{4}{3}$ del diapente e del diatessaron --

$$
\begin{array}{l}
\hspace{4.2truecm} \overbrace{\hspace{2.2truecm}}^{epogdoon} \hspace{3truecm}	\\
\underbrace{\underbrace{	1\;({\rm VI})\qquad \dfrac{9}{8}
		\qquad 
		\dfrac{81}{64} \qquad \dfrac{4}{3}\;({\rm VIII})}_{diatessaron} 
	\qquad \dfrac{3}{2}\;({\rm VIIII})}_{diapente}
		\qquad \dfrac{27}{16} \qquad \dfrac{243}{128}\qquad 2\;({\rm XII})\\
\hspace{4truecm}\underbrace{\hspace{1.9truecm}\underbrace{\hspace{5.6truecm}}_{diatessaron}}_{diapente}\\
\underbrace{\hspace{11.5truecm}}_{diapason}\\
\underbrace{\hspace{1.7truecm}}_{epogdoon} \;
\underbrace{\hspace{1.2truecm}}_{epogdoon} \;
\underbrace{\hspace{1truecm}}_{256:243}
 \hspace{2.2truecm}
\underbrace{\hspace{1.8truecm}}_{epogdoon}\;
\underbrace{\hspace{1.4truecm}}_{epogdoon}\;
\underbrace{\hspace{1.5truecm}}_{256:243} 
\end{array}
$$
(i numeri romani rimandano al suono analogo di origine). In base alla numerazione semplificata, 
il fondamento dei quattro suoni consonanti dichiarati sulla tavoletta corrisponde all'insieme
\begin{equation}
	\label{tav}
	{\cal T}=\left\{1, \dfrac{4}{3}, \dfrac{3}{2}, 2\right\}.
\end{equation}
Lo schema, pur privo del fascino cinquecentesco, riproduce gli intervalli della tavoletta ed aggiunge quelli che intercorrono fra suoni consecutivi, disposti sull'ultima riga, per dare origine alla {\it scala pitagorica}
\begin{equation}
\label{scalapit}
{\mathfrak S}_{\cal P}=\left\{1, \dfrac{9}{8}, \dfrac{81}{64}, \dfrac{4}{3}, \dfrac{3}{2}, \dfrac{27}{16}, \dfrac{243}{128}, 2\right\}.
\end{equation}

\noindent
Se pu\`o far comodo una collocazione sonora pi\`u moderna, la sequenza riproduce la probabilmente familiare scala diatonica DO--RE--MI--FA--SOL--LA--SI--DO, nell'intonazione pitagorica.
I commenti sono molteplici: il diatessaron abbraccia effettivamente quattro suoni (estremi compresi) ed il diapente cinque, la mutua distanza \`e prevalentemente l'epogdoon, che tuttavia non si incastra perfettamente nel diapason; esordisce un intervallo pi\`u piccolo di misura $\frac{256}{243}$, che conserva un'impronta della {\it tetraktys} pitagorica nella scrittura $\frac{2^8}{3^5}$ ma mostra al tempo stesso una struttura complicata, se pensiamo alla suddivisione della corda. Questo intervallo minimo \`e il cosiddetto {\it limma} (ovvero ``residuo'') ed interviene nell'articolata questione, qui estranea, di estendere la scala a 13 suoni (estremi compresi), suggerita in modo forzato ed antistorico dall'odierna scala cromatica, quella che si ascolta suonando in progressione i tasti (bianchi e neri) di un pianoforte
\footnote{Dal punto di vista della matematica, la questione va associata alla non finitezza della successione di suoni che si generano procedendo per diapente: la necessit\`a di una selezione finita di suoni ha prodotto svariate soluzioni e compromessi, soprattutto in sede di accordatura di strumenti.
Si pu\`o per esempio verificare (vedi \cite{tal}) che 12 diapente in avanti e 12 all'indietro stabiliscono una situazione di quasi chiusura e completamento, nel senso che al dodicesimo passo in avanti si trova un suono molto vicino a quello di partenza, al dodicesimo passo all'indietro si trova un suono molto vicino al diapason; inoltre, i 26 suoni ottenuti si dispongono a coppie per fornire una buona approssimazione per difetto e per eccesso di ciascuno dei 13 suoni della scala cromatica del pianoforte:
invertendo il punto di vista, \`e come se ciascun suono di tale scala attuasse una sintesi dei due esemplari pitagorici calante e crescente che lo abbracciano.}.
La sequenza degli otto suoni poco sopra esibita \`e gi\`a sufficiente per proseguire l'argomento. 

\subsection{I suoni naturali}

\noindent
Accogliamo una novit\`a: valichiamo il precetto pitagorico di intervento dei soli numeri 1, 2, 3 e 4 nell'armonia e proclamiamo la consonanza del suono prodotto dalla corda 
divisa in cinque parti. Entrano a far parte delle consonanze (sorvolando su complicate gerarchie) i suoni prodotti da $1/5$ della corda, $2/5$, $3/5$ e $4/5$.
Il naturale e disarmato passaggio da 4 a 5 si articola in realt\`a attraverso un lungo e serio dibattito di epoca tardomedievale, quando i profondi mutamenti della pratica musicale 
ponevano dubbi sull'uso di alcuni suoni pitagorici\footnote{
L'assestamento teorico dell'innovazione avviene attraverso i trattati di Gioseffo Zarlino (\cite{zarlino})
nel XVI secolo, nell'epoca di forte convinzione sull'esistenza di rapporti universalmente validi, svelati dalla musica e adoperati dall'architettura. \`E proprio la ripetuta presenza di proporzioni nella cosiddetta {\it scala naturale} promossa da Zarlino che verr\`a messa in evidenza pi\`u avanti.}.
I suoni inediti sono in realt\`a solamente due: quello prodotto dalla lunghezza $4/5$ (oppure $1/5$ o $2/5$ a meno di diapason) e quello prodotto dalla lunghezza $3/5$; i corrispondenti suoni hanno altezza, nel diapason normalizzato, $5/4$ e $5/3$. Cos\`\i$\,$come la divisione fino a quattro parti della corda 
d\`a luogo, al netto di distanze lunghe un diapason, ai quattro suoni (\ref{tavoletta}), la divisione fino a cinque parti aggiorna l'insieme sonoro di base a
\begin{equation}
	\label{div5}
{\cal T}_5=\left\{ 1, \frac{5}{4}, \frac{4}{3}, \frac{3}{2}, \frac{5}{3}, 2  \right\} 
\end{equation}
a meno delle medesime identificazioni.
L'ampliamento dei suoni ${\cal T}_5$ segue un percorso ben differente da quello dei suoni pitagorici, come si evince dai manuali di acustica musicale o da interventi in merito presenti in letteratura: 
non c'\`e da aspettarsi un procedimento ciclico, come per la scala pitagorica, per generare nuovi suoni;
in effetti le molteplici possibilit\`a che si aggiungono con l'apertura alla divisione per cinque rendono complicato questo approccio.  
Il modo di procedere si basa piuttosto su un aggiustamento -- o per meglio dire un {\it temperamento} --  di alcuni tra i suoni pitagorici.
Il nucleo originale (\ref{tav}) delle consonanze non viene comprensibilmente toccato, non solo per rispetto ``ideologico'', ma anche per la non necessit\`a di alterare questi semplici rapporti numerici.
Diversamente, la vicinanza di $\frac{81}{64}$ e $\frac{27}{16}$
a $\frac{5}{4}$ e $\frac{5}{3}$ convince ad aggiornare le due complicate frazioni pitagoriche mediante questi ultimi due suoni a disposizione. L'altro elemento con medesima sorte \`e il suono pitagorico $\frac{243}{128}$, prossimo a $\frac{15}{8}$, 
che \`e utilizzabile grazie al nuovo assortimento procurato dal cinque e che sostituisce il primo nella nuova sequenza. Si stabilisce la scala dei cosiddetti {\it suoni naturali} con il riconoscere i valori (fissati nuovamente e comodamente nel diapason tra 1 e 2)
\begin{equation}
\label{scalanat}
{\mathfrak S}_{\cal N}=\left\{1, \dfrac{9}{8}, \dfrac{5}{4}, \dfrac{4}{3}, \dfrac{3}{2}, \dfrac{5}{3}, \dfrac{15}{8}, 2\right\}.
\end{equation}
che dispone ed amplia l'insieme (\ref{div5}) in accordo al seguente diagramma, in cui 
continuiamo a riportare anche l'antica presenza degli intervalli della tavoletta:
$$
\begin{array}{l}
	\hspace{4.4truecm} \overbrace{\hspace{1.8truecm}}^{epogdoon} \hspace{3truecm}	\\
	\underbrace{\underbrace{	1\hspace{1.4truecm}\dfrac{9}{8}
		\hspace{1.2truecm}
			\dfrac{5}{4} \hspace{1.truecm} \dfrac{4}{3}}_{diatessaron} 
	 \hspace{1.4truecm}  \dfrac{3}{2}}_{diapente}\hspace{1.2truecm} \dfrac{5}{3}
	\hspace{1.4truecm} \dfrac{15}{8} \hspace{1.truecm}  2\\
	\hspace{4.6truecm}
	\underbrace{\hspace{1.5truecm}\underbrace{\hspace{4.7truecm}}_{diatessaron}}_{diapente}\\
	\underbrace{\hspace{10.8truecm}}_{diapason}\\
	\underbrace{\hspace{1.7truecm}}_{t.~mag.} \;
	\underbrace{\hspace{1.4truecm}}_{t.~min.} \;
	\underbrace{\hspace{1.1truecm}}_{s.~mag.}\;
	\underbrace{\hspace{1.6truecm}}_{t.~mag.}\;
	\underbrace{\hspace{1.35truecm}}_{t.~min.}\;
	\underbrace{\hspace{1.65truecm}}_{t.~mag.}\;
	\underbrace{\hspace{1.1truecm}}_{s.~mag.} 
\end{array}
$$

\vspace{.7truecm}

\noindent
La denominazione attuale dei suoni (\ref{scalapit}) \`e la medesima data in precedenza, ovvero le note DO--RE--MI--FA--SOL--LA--SI--DO, nell'intonazione stavolta naturale. Va tuttavia osservato -- e su ci\`o torneremo nelle conclusioni -- che, pur essendo antica la denominazione richiamata\footnote{Guido monaco, attorno all'anno Mille, utilizz\`o sei suoni in progressione dell'inno ``Ut queant laxis'' di Paolo Diacono, di epoca longobarda, per nominare gli intervalli dell'{\it esacordo musicale}
UT(=DO)--RE--MI--FA--SOL--LA; il SI \`e assente.}, alla base del ragionamento della teoria musicale cinquecentesca -- e a maggior ragione del mondo antico -- c'\`e 
l'intervallo, ovvero il mutuo rapporto fra i suoni -- gi\`a noto dall'antichit\`a prima ancora di dare il nome alle note -- e non la scala, la sequenza sonora in successione crescente.
Quest'ultima, supporto inevitabile attuale per trattare i presenti argomenti, porta con s\'e la moderna impronta sonora della {\it tonalit\`a}, dell'evoluto clima musicale che fa da contorno e detta le regole alle note, faccenda che non interviene nel nostro percorso.

\noindent
Tornando a (\ref{scalanat}), il commento \`e indubbiamente quello di una semplificazione delle frazioni che fissano i suoni e delle misure degli intervalli: giustappunto, la serie riceve la denominazione di {\it 
suoni naturali}\footnote{La spontanea analogia con il fenomeno dei suoni naturali che accompagnano il suono fondamentale richiamato all'inizio, pur comportando parziali concidenze, \`e antistorica.}. \`E aumentata la gamma delle distanze, fatto non necessariamente sfavorevole: oltre all'epogdoon (che prende il nome di {\it tono maggiore}, abbreviato con {\it t.~mag.~}) di misura $9/8$, si trova la distanza $10/9$ ({\it tono minore}, {\it t.~min.~}) e quella minima $16/15$ ({\it semitono maggiore}, {\it s.~mag.~}) che rimpiazza e semplifica il limma. 

\noindent
Altri elementi di riflessione circa la formazione della scala lasciano pi\`u incertezze: da una parte, il criterio ``manuale'' di rimpiazzare i suoni appare lontano da un sostegno matematico; dall'altra, non si vede la ragione musicale per cui un suono cifrato da una frazione pi\`u semplice sia preferibile ad un altro, se non per l'agevolazione di riprodurli suddividendo la corda. Gli interventi della matematica a sostegno della selezione dei suoni naturali saranno oggetto della Sezione 3.

\noindent
Ribaltando l'ordine di apparizione precedente, mostriamo solo ora come Zarlino ha completato il diapason, in una celebre illustrazione del Trattato \cite{zarlino}, nella quale si riconoscono, in mezzo ad uno schema pi\`u completo, le indicazioni fornite.

\begin{figure}
	\includegraphics[width=0.36\textwidth]{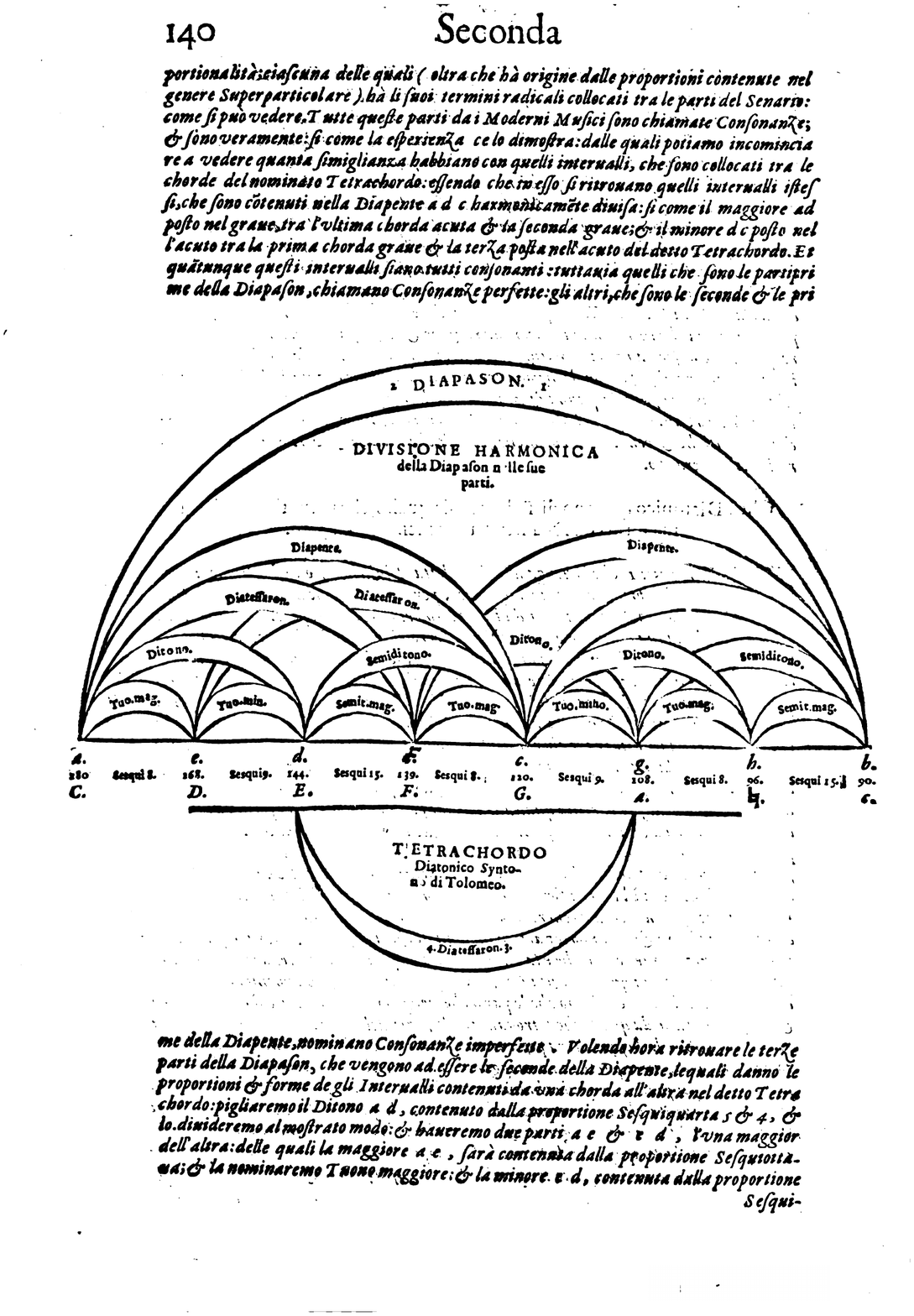}
	\caption{Prospetto della suddivione del diapason, secondo la teoria dei suoni naturali di Zarlino. La posizione del suono \`e indicata mediante la lunghezza della corda, decrescente verso l'acito. L'immagine \`e tratta da \cite{zarlino}.}
\end{figure}

\subsection{I suoni equabili}

\noindent
Innanzi al diapason, alla nozione di distanza e alla regola moltiplicativa per gli intervalli, il matematico non pu\`o trattenersi dall'affermare la soluzione pi\`u logica, che consiste nel disporre i suoni {\it in modo equidistante}, ovvero suoni separati da intervalli di medesima lunghezza. L'impostazione matematica propone la ricerca dei coefficienti  
$1=\alpha_1<\alpha_2<\dots<\alpha_N<2$ che distribuiscano i suoni $\alpha_1\nu=\nu$, $\alpha_2\nu$, $\dots$, $\alpha_N\nu$, $2\nu$ -- appartenenti al diapason tra $\nu$ e $2\nu$ -- alla medesima distanza. La regola moltiplicativa della distanza
comporta
$$
\dfrac{\alpha_2}{1}=\dfrac{\alpha_3}{\alpha_2}=\;\dots\;=\dfrac{\alpha_N}{\alpha_{N-1}}=\dfrac{2}{\alpha_N}
$$
che dispone gli $N$ numeri $\alpha_1=1$, $\alpha_2$, $\dots$, $\alpha_N$ in progressione geometrica di ragione $\sqrt[N]{2}$ e con elementi pari a 
\begin{equation}
	\label{scalatemp}
	\alpha_\kappa=\sqrt[N]{2^{\kappa-1}}, \qquad \kappa=1, \dots, N.
\end{equation}
A dire il vero, la leggera presentazione, a mo' di esercizio, dei suoni equidistanziati \`e irrispettosa verso la pesante problematica di uniformare i suoni per fini pratici, a costo di un ipotetico appiattimento sonoro causato dall'abbattimento delle variet\`a di intervalli.

\noindent
D'altra parte, la spontanea domanda di suoni equidistanti trova in realt\`a risposta solo nella teoria, dal momento che i numeri (\ref{scalatemp}) sono irrazionali ed il preciso valore matematico non \`e nella pratica eseguibile, se non accogliendo un certo grado di approssimazione; dalla lista (\ref{scalatemp}) 
spariscono il diapente ed il diatessaron. Per includere nella scala equiripartita suoni prossimi a quelli pitagorici e naturali dapprima discussi, va scelto  
$N=12$ per trovare la cosiddetta {\it scala del temperamento equabile} formata da 12 unit\`a ({\it semitoni}) uguali ed il cui principale vantaggio consiste nell'essere identica, da qualunque suono di partenza $\nu$ si proceda\footnote{Chiamando senza troppo impegno la nota $\nu$  {\it tonalit\`a}, il vantaggio pratico dell'invariabilit\`a risiede nel poter comporre in ogni totalit\`a e nel passare pi\`u agevolmente da una tonalit\`a ad un'altra. I dodici suoni dell'ottava di un pianoforte (sette tasti bianchi e cinque neri) fanno riferimento alla scala temperata ed in fase di accordatura l'avvicinamento ai suoni teorici avviene sfruttando il fenomeno dei battimenti oppure utilizzando apparecchi elettronici. Diversamente,  gli strumenti a tastiera di origine pi\`u antica (clavicembalo, clavicordo, organo, ...) frequentemente non fanno uso del temperamento equabile.}.

\noindent
Commenti e confronti fra le tre tipologie di suoni procurate -- pitagorici, naturali ed equabili -- non fanno parte dei propositi prefissi e sono comunque argomenti con vasto repertorio, anche divulgativo. 
\`E parimenti omesso per le medesime ragioni un accenno alla pur interessante storia del temperamento equabile, a partire dall'immancabile precursore nel mondo antico -- Aristosseno di Taranto\footnote{Nell'affresco della Scuola di Atene Aristosseno viene identificato con il filosofo in piedi di fronte a Pitagora, oppure con quello seduto dietro il medesimo; in qualsiasi modo, entrambi sono in atto di fissare le azioni di Pitagora e di prendere appunti su un libro; l'interpretazione \`e probabilmente guidata dalla contrapposizione delle rispettive teorie musicali.} --, attraverso il sostegno di Vincenzo Galilei in polemica con Zarlino, l'attenzione di matematici come Stevino e Cartesio, fino all'ingegnosa sistemazione teorica e pratica di Andreas Werckmeister, alla fine del 600.

\section{I medi proporzionali}

\noindent
Inquadriamo ora un aspetto della tavoletta della Figura 1 che fa da preludio alle successive considerazioni: 
\`e possibile rintracciare in essa la presenza schematica dei tre medi proporzionali pi\`u familiari in matematica, definiti come segue, con il pratico supporto di due segmenti di lunghezza $a$ e $b \geq a$.

\noindent
Nella {\it proporzione aritmetica} un terzo segmento di lunghezza $m_A$ si inframezza in modo che superi il primo segmento tanto quanto viene superato dal secondo:
$$
m_A-a=b-m_A
$$
Nella {\it proporzione geometrica} il terzo segmento di misura $m_G$ supera il primo formando un rapporto con esso pari a quello che esso forma con il secondo segmento:
$$
a:m_G =m_G:b
$$
Infine, nella {\it proporzione armonica} il segmento di misura $m_H$ si interpone fra $a$ e $c$ in modo da superare il primo ed essere superato dal secondo rispettando il medesimo rapporto dei due segmenti assegnati:

$$
\dfrac{m_H-a}{b-m_H}=\dfrac{a}{b}
$$
Questa presentazione dei medi proporzionali, dal sapore antico, gi\`a appare nei Commenti che Marsilio Ficino aggiunse alla traduzione latina del dialogo platonico {\it Timeo} (\cite{timeo}). La predominanza dei medi proporzionali e la corrispondenza fra suoni ed elementi architettonici che pervade il periodo rinascimentale, dall'Alberti a Palladio, trae una delle principali fonti di ispirazione per l'appunto dal testo del Ficino.
Proprio nel {\it Timeo}\footnote{Di nuovo un richiamo all'affresco di Raffaello: l'attempato filosofo Platone--Leonardo da Vinci al centro della scena stringe fra il braccio ed il fianco il testo del Timeo.} viene offerta una spiegazione sul fatto che i tre medi proporzionali individuano tutti gli intervalli della scala musicale: sostanzialmente, la motivazione sviluppa  la stilizzata comparsa dei medi nella tavoletta, come dichiarato dapprima, che consiste nella proporzione geometrica di cui fanno parte i diapason ${\rm VI}$ e ${\rm XII}$ (a rappresentare i raddoppi e i dimezzamenti), nel medio aritmetico ${\rm VIIII}$ degli estremi ${\rm VI}$ e ${\rm XII}$, nel medio armonico ${\rm VIII}$ dei medesimi estremi.

\noindent
Ecco il passo in avanti da compiere, sommessamente evocato dalla tavoletta: non solo cercare l'origine, la motivazione numerica del singolo suono, ma raccordare con un elemento intermedio
i suoni tra loro, trovare le motivazioni dell'armonia nella presenza di medi proporzionali.

\noindent
E questo aggiorna il motivo principale dello scritto che portiamo avanti: il suono \`e reso intellegibile dal numero, dal rapporto che lo pone in relazione con il tutto. Ma ciascuna parte del tutto deve accordarsi con l'altra, ed il mezzo primario di indagine risiede nella matematica delle concatenazioni proporzionali.

\noindent
Se nell'Alberti c'\`e la piena coscienza ed affermazione dell'identificazione tra rapporto spaziale e musicale, \`e con i rapporti palladiani\footnote{Il trattato {\it ``I quattro libri di Architettura''}
del 1570 \`e ritenuto il punto di riferimento pi\`u elevato per un corretto sistema costruttivo proporzionale.} che l'architettura e la teoria musicale trovano un'inedita ed irripetibile intesa. I moduli architettonici vengono concatenati come accordi, le proporzioni dei suoni e dello spazio obbediscono ad un medesimo ed unico sistema armonico, a rapporti universalmente validi che la musica ha il privilegio di esibire. Palladio sceglie ambienti e misure come se fossero note; dispone, organizza, raccorda gli ambienti come se formasse una scala musicale.

\noindent
Il principio costruttivo ``sinfonico'' -- nel senso etimologico di ``suono insieme'' -- e pi\`u in generale il nesso tra l'architettura rinascimentale e i coevi elementi culturali sono mirabilmente 
narrati con disciplina e fascino nel celebre testo \cite{witt} di Rudolf Wittkower\footnote{Berlino 1906 -- New York 1971, celebre storico dell'architettura e dell'arte, ha chiarito l'evidenza del modello rinascimentale delle proporzioni, contrariamente ad una interpretazione puramente estetica.}: in particolare, nella Parte quarta ``Il problema della proporzione armonica in architettura'' l'Autore del saggio ci accompagna in celebri ville palladiane a spasso per i vari ambienti, decifrandone con rigore le varie proporzioni e gli scambievoli rapporti armonici. Quell'apparentemente indecifrabile percezione di ordine, simmetria, bellezza -- insomma di benessere per la mente -- destata dall'ammirazione di ambienti del genere come pure dall'ascolto di suoni armoniosi, incontra la logica spiegazione, la naturale delucidazione -- la scienza -- nelle insuperate argomentazioni del Wittkower.
L'armonico progetto delle stanze della superba Villa Barbaro\footnote{Per l'umanista Daniele Barbaro, autore di importanti commenti sull'architettura di Vitruvio, Andrea Palladio costru\`i$\,$ a Maser intorno al 1558 la magnifica villa di Maser, decorata dagli splendidi affreschi di Paolo Veronese.} riportato in figura 3 ed il prospetto degli intervalli della scala naturale della figura 2
nella loro essenza sono identici, ispirati dal medesimo principio universale.
\begin{figure}
	\includegraphics[width=0.36\textwidth]{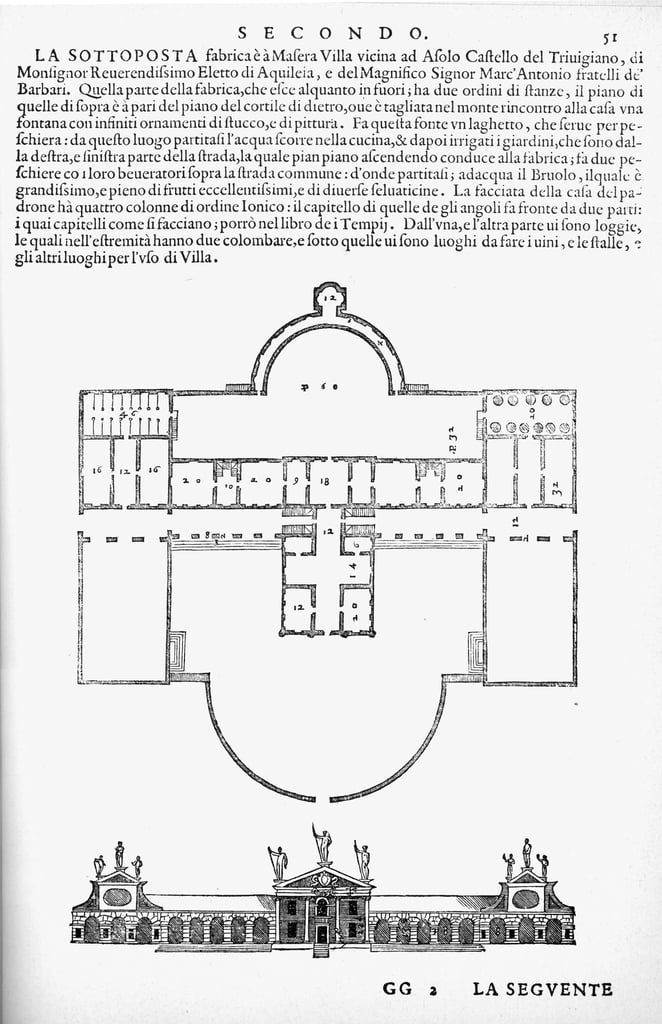}
	\caption{Planimetria e facciata della Villa Maser. La pagina \`e tratta dalla celebre opera ``I quattro libri dell'Architettura'' di Andrea Palladio, 1570.}
\end{figure}
\begin{figure}
	\includegraphics[width=0.36\textwidth, angle=270]{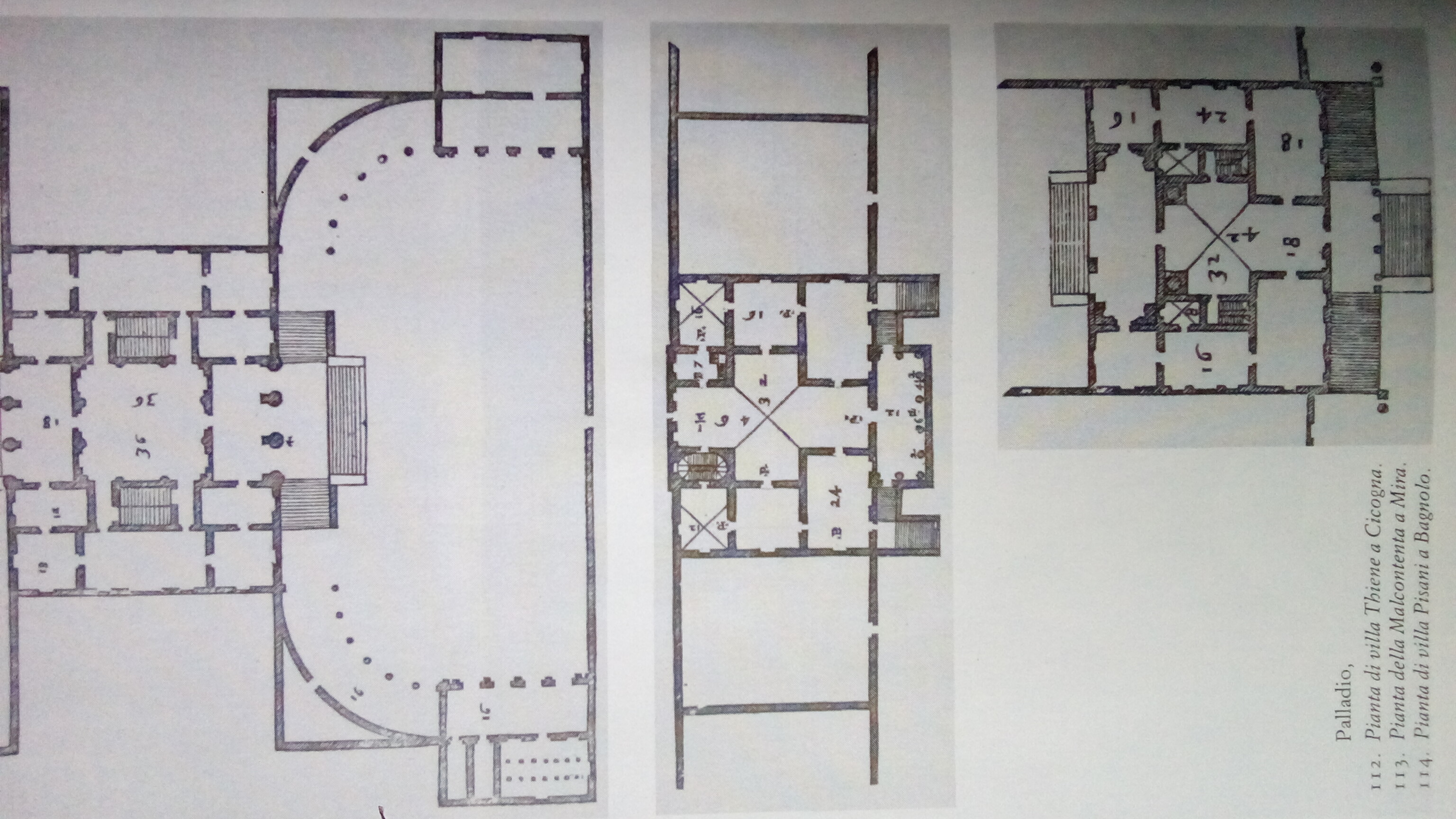}
	\caption{Altre piante di ville del Palladio, concepite come universi musicali. L'immagine \`e tratta da \cite{witt}.}
\end{figure}

\noindent
Sul terreno dei princip{\^\i} universali che gettino la luce dell'assoluto su qualche attivit\`a umana si aggira perennemente -- e pi\`u o meno proficuamente -- il matematico: 
cambiamo decisamente scenario e procediamo da matematici, per compiere un'indagine neutrale, razionale, cieca (anzi sorda) al cospetto di motivazioni musicali che potrebbero deviarne il percorso.
Le formule esplicite, ricavabili algebricamente dalle definizioni dei medi geometrici,
\begin{equation}
	\label{mediaag}
	m_A(a,b)=\dfrac{a+b}{2}, \quad m_G(a,b) = \sqrt{a b}
\end{equation}

\begin{equation}
	\label{mediah}
	m_H(a,b) = \dfrac{2ab}{a+b}=\dfrac{2}{\frac{1}{a}+\frac{1}{b}}
\end{equation}
danno agevole accesso ad alcune immediate propriet\`a:
\begin{description}
	\item[$\bullet$] le tre medie rispettano il rapporto di similitudine, nel senso che le medie dei segmenti di lunghezza $\lambda  a$ e $\lambda b$ sono le (\ref{mediaag}) e (\ref{mediah}) moltiplicate per $\lambda$, qualunque sia il numero reale positivo $\lambda$;
	\item[$\bullet$]  $m_H\leq m_G \leq m_M$ e l'uguaglianza vale, in entrambi i segni, se e solo se $a=b$, caso nel quale le tre medie coincidono con $a$;
	\item[$\bullet$] la prima uguaglianza in (\ref{mediah}) equivale alla relazione $m_A m_H = m_G^2$, ovvero ciascuna media \`e ricavabile dalle altre due;	in particolare $m_H=\frac{m_G^2}{m_A}$; 
	\item[$\bullet$] la seconda uguaglianza in (\ref{mediah}) mostra che la media armonica di $a$ e $b$ coincide con il reciproco della media aritmetica dei reciproci $1/a$ e $1/b$:
	$$
	m_H(a,b)=m_A^{-1}\left(\frac{1}{a}, \frac{1}{b}\right), \qquad m_A(a,b)=m_H^{-1}\left(\frac{1}{a}, \frac{1}{b}\right);
	$$
	\item[$\bullet$] i reciproci di una proporzione geometrica formano una proporzione geometrica e 
	$$m_G\left( \frac{1}{a}, \frac{1}{b} \right)=\frac{1}{m_G(a,b)}.
	$$
\end{description}

\noindent
Se i segmenti sono corde poste in vibrazione, la proporzionalit\`a inversa tra la lunghezza $\ell$ della corda e la frequenza di vibrazione $\nu$
\begin{equation}
	\label{fil}
	\nu=\dfrac{\kappa}{\ell}
\end{equation}
(la costante $\kappa$ \`e la medesima a parti\`a di densit\`a e di tensione) stabilisce la seguente
\begin{propr}
Se due corde di lunghezza $a$ e $b$ hanno media armonica $m_H$, allora la frequenza di vibrazione della corda lunga $m_H$ \`e la media aritmetica delle frequenze delle due corde.
\end{propr}

\noindent
{\bf Dim.} Le frequenze di vibrazione, in virt\`u della (\ref{fil}), sono
$\nu_a=\dfrac{\kappa}{a}$, $\nu_b =\dfrac{\kappa}{b}$; d'altra parte  dalla (\ref{mediah}) si deduce
$$
m_H = \dfrac{2\kappa}{\frac{\kappa}{a}+\frac{\kappa}{b}}=\dfrac{2\kappa}{\nu_a +\nu_b}= \dfrac{\kappa}{m_A(\nu_a, \nu_b)}
$$
dunque  $\nu_{m_H}= \dfrac{\kappa}{m_H}=m_A (\nu_a, \nu_b)$. $\quad \square$

\subsection{Alla ricerca di proporzioni}

\noindent
L'insegnamento da scorgere nell'ispirazione ai canoni dell'architettura rinascimentale \`e la predilezione per la proporzione armonica delle dimensioni spaziali, che giocano il ruolo di corde musicali: 
la prima esplorazione da compiere \`e riscontrare la presenza di terne di corde formanti proporzioni armoniche
nei sistemi di suoni dapprima delineati.
La Propriet\`a 1 fa trasferire l'indagine sulle proporzioni armoniche 
alla ricerca di proporzioni aritmetiche sulle frequenze o su equivalenti quantit\`a inversamente proporzionali alle lunghezze delle corde\footnote{Questo fatto, non in termini di frequenza ma di inverso della lunghezza, \`e presente gi\`a nel Trattato \cite{zarlino}, in cui Zarlino mostra che le consonanze siano determinate ``{\it dal medio tanto aritmetico che armonico}''.}.
Richiamiamo le due sequenze (\ref{scalapit}) e (\ref{scalanat}) di suoni pitagorici e suoni naturali\footnote{La circostanza dei suoni equabili verr\`a commentata pi\`u in seguito.} dapprima delineate, per osservare che dal punto di vista del numero esse appartengono all'universo dei suoni prescritto dalle divisioni consentite:
$$
\begin{array}{ll}
{\mathfrak S}_{\cal P}\subset {\cal P}, & {\mathfrak S}_{\cal N}\subset {\cal N}
\end{array}
$$
dove
\begin{equation}
	\label{suonipit}
{\cal P}=\{2^m 3^n, \;\; m,n \in {\Bbb Z} \}
\end{equation}

\begin{equation}
	\label{suoninat}
{\cal N}=\{	2^m 3^n 5^p, \;\; m,n,p \in {\Bbb Z} \}
\end{equation}
con ${\Bbb Z}$ l'insieme dei numeri interi.
L'illustrazione delle medie aritmetiche tra gli elementi di ${\mathfrak S}_{\cal P}$ o di ${\mathfrak S}_{\cal N}$, esibite nelle due tabelle che seguono (il calcolo proviene in ciascuna posizione dal corrispondente incrocio degli elementi incontrati sulla prima riga e prima colonna),

\vspace{.5truecm}
 
 \begin{tabular}{c||c|c|c|c|c|c|c|c|}
${\mathfrak S}_{\cal P}$	& 1	& 9/8 &  81/64& 4/3& 3/2& 27/16& 243/128 & 2\\
 	\hline \hline 
 1	&	& 17/16 &145/128 &7/6 &\cellcolor[gray]{0.5} 4/3& 43/32& 371/128& \cellcolor[gray]{0.5}3/2\\
\cline{1-1} 	\cline{3-9}
9/8 	&	&  & 153/128& 59/48& 21/16& 45/32& 387/256& 25/16\\
 \cline{1-1} 	\cline{4-9}
81/64 	&	&  & & 499/384&177/256 & 189/128& 405/256& 209/128\\
\cline{1-1}  	\cline{5-9}
4/3 	&	&  & & & 17/12& 145/96& 1241/384& 5/3\\
\cline{1-1}  	\cline{6-9}
3/2 	&	&  & & & & 51/32& 435/256& 7/4\\
\cline{1-1}  	\cline{7-9}
27/16 	&	&  & & & & & 459/256& 59/32\\
\cline{1-1}  	\cline{8-9}
243/128 	&	&  & & & & & & 499/256\\
 	\hline
 \end{tabular}

\vspace{.5truecm}

\begin{tabular}{c||c|c|c|c|c|c|c|c|}
	${\mathfrak S}_{\cal N}$ & 1	& 9/8 & 5/4& 4/3& 3/2& 5/3& 15/8 & 2\\
	\hline \hline 
	1 	&	& 17/16 &\cellcolor[gray]{0.5}9/8 & 7/6&\cellcolor[gray]{0.5}5/4 & \cellcolor[gray]{0.5}4/3&23/16 & \cellcolor[gray]{0.5}3/2\\
	\cline{1-1} 		\cline{3-9}
	9/8	&	&  & 19/16& 59/48& 21/16& 67/48& \cellcolor[gray]{0.5}3/2&\cellcolor[gray]{0.9}25/16 \\
	\cline{1-1} 			\cline{4-9}
	5/4	&	&  & &31/24 & 11/8& 35/24& \cellcolor[gray]{0.9}25/16& 13/8 \\
	\cline{1-1} 			\cline{5-9}
	4/3	&	&  & & &  17/12&\cellcolor[gray]{0.5}3/2& 77/48&\cellcolor[gray]{0.5}5/3 \\
	\cline{1-1} 		\cline{6-9}
	3/2	&	&  & & & &19/12 &\cellcolor[gray]{0.9}27/16& 7/4\\
	\cline{1-1} 			\cline{7-9}
	5/3	&	&  & & & & & 85/48& 11/6\\
	\cline{1-1} 		\cline{8-9}
	15/8	&	&  & & & & & &31/16 \\
	\hline
\end{tabular}

\vspace{.5truecm}

\noindent
\`e gi\`a di per s\'e espressiva. Con il fondo pi\`u scuro si sono evidenziate le medie che non escono dalla scala  ${\mathfrak S}_{\cal P}$ nella prima tabella, o da ${\mathfrak S}_{\cal N}$ nella seconda,  con il fondo grigio meno scuro le medie estranee alla scala, tuttavia appartenenti a ${\cal P}$ (ovvero di tipo (\ref{suonipit})) nel primo caso, appartenenti a ${\cal N}$ (dunque di tipo (\ref{suoninat})) nel secondo. Nella prima tabella l'unica presenza di proporzioni aritmetiche \`e limitata alla stretta cerchia dei quattro suoni (\ref{tavoletta}), con le due proporzioni evidenziate; il resto delle caselle presenta dei rapporti talvolta complicati e non connaturati con la corrispondente ricerca di  ``moduli costruttivi'' da sezioni della corda.

\noindent
Ben differente appare la situazione nella seconda tabella: ogni valore della scala ${\mathfrak S}_{\cal N}$ \`e interessato da almeno una proporzione, inoltre vengono generati due suoni naturali (fondo grigio chiaro), ovvero del tipo (\ref{suoninat}).

\noindent
La predisposizione dei suoni naturali alla proporzione armonica (intesa sempre nel senso della Propriet\`a 1) emersa dalla tabella \`e avvalorata anche da pi\`u di un metodo costruttivo che possiamo architettare per giungere ai suoni ${\mathfrak S}_{\cal N}$ procedendo per proporzioni, senza che sia neppure necessario il ``dogma'' (\ref{tavoletta}): pu\`o valere la pena di 
esporre il seguente, che presenta dei tratti quasi esclusivamente matematici. 
A partire dal solo suono di base, corrispondente a $1$, e dal suo diapason $2$, si aggiungono in successione le medie 
\begin{equation}
\label{riemp}
m_A(1,2)=\frac{3}{2}\;\;\rightarrow\;\;m_A(1, \frac{3}{2})=\frac{5}{4}\;\;\rightarrow\;\;m_A(1, \frac{5}{4})=\frac{9}{8}
\end{equation} legittimamente appartenenti all'universo dei suoni (\ref{suoninat}); per contro, il frazionamento successivo $m_A(1,\frac{9}{8})=\frac{17}{16}$ farebbe uscire dai suoni naturali. 
Il temporaneo insieme di suoni  $\{ 1, \frac{9}{8}, \frac{5}{4},  \frac{3}{2}, 2\}$ pu\`o essere ampliato sempre ricorrendo alle proporzioni, cercando di compensare alcune distanze eccessive tra suoni consecutivi:
\begin{description}
	\item[$(i)$]
l'unico suono naturale compreso fra $\frac{3}{2}$ e $2$ in grado di formare una proporzione aritmetica con i quattro precedenti \`e  $\frac{15}{8}$, tramite l'espressione $m_A(\frac{5}{4}, \frac{15}{8})=\frac{3}{2}$. 
\item[$(ii)$] 
una seconda considerazione, ancora a sostegno di una sorta di unicit\`a della costruzione, \`e la seguente: i suoni naturali che possono formare medie aritmetiche con gli estremi del diapason, nel senso di due suoni
$\sigma_1$ e $\sigma_2$ per cui $\sigma_1= m_A(1,\sigma_2)$ e $\sigma_2=m_A(\sigma_1, 2)$, sono solamente $\frac{5}{4}$ e $\frac{5}{3}$.
\end{description}
Intervallare il diapason con i tre suoni (\ref{riemp}) e con i tre emersi da $(i)$ e $(ii)$ 
fa conseguire esattamente i suoni naturali ${\mathfrak S}_{\cal N}$; tuttavia la soddisfazione \`e attenuata da pi\`u di una perplessit\`a: innanzitutto, la decisione di fermarsi \`e 
principalmente motivata dall'aver ottenuto qualcosa che, in un certo senso, non ha tradito le aspettative, ovvero le sette note. La disinteressata neutralit\`a della matematica, ignara della teoria musicale, non vedrebbe la ragione per cui non procedere al calcolo di ulteriori medie ``naturali'', effettivamente esistenti. 
In secondo luogo, un oggetto di critica -- in cui rientra peraltro l'obiezione appena esposta -- \`e quello di aver seguito un metodo non sistematico, la cui consequenzialit\`a si aggiusta pi\`u su considerazioni sensate ma ``esterne'', anzich\'e su una linea precisa, anteposta.
 
\subsection{Un metodo che genera medi proporzionali}

\noindent
A questo punto, anzich\'e cercare un procedimento 
che insegua la matematica armonia {\it in itinere} -- laddove non mancano  considerevoli esempi di illustri esponenti del pensiero teorico, quali Mersenne (\cite{mersenne}), Cartesio (\cite{cartesio}), Rameau (\cite{rameau}), Eulero (\cite{eulero}) -- meditiamo sulla capacit\`a dell'armonica matematica di portare un {\it incipit}, un avvio, ad un {\it exitus} convincente, un prodotto musicale, attraverso la conduttura rigorosa di un algoritmo.
Rientriamo cos\`\i$\,$ nell'esatto perimetro della matematica per porre la seguente questione: \`e possibile formulare un procedimento sulla linea intenzionale di accostare i suoni sulla base delle proporzioni, ma che proceda senza intromissioni, per offrire alla fine un insieme di suoni ``ideale'' 
in merito al reciproco accordo che le proporzioni medesime trasmettono?

\noindent
Per fissare un'idea dal punto di vista matematico, 
procurandoci uno strumento per conteggiare
la presenza e la generazione di medie aritmetiche, definiamo il {\it generatore di medie}  
\begin{equation}
	\label{genmedie}
{\Bbb M}\left\langle \; {\cal I} \right \rangle_{\mathfrak R} \qquad {\cal I}\;\;\textrm{insieme di numeri}, \quad {\mathfrak R}\;\;\textrm{restrizioni}
\end{equation} 
l'insieme di tutte le possibili medie aritmetiche calcolate con le coppie di elementi in ${\cal I}$, accettando il risultato solo se \`e soddisfatta la restrizione ${\mathfrak R}$. 
La destinazione \`e prevedibile: ${\cal I}$ \`e un complesso originario di suoni, 
${\mathfrak R}$ \`e la richiesta che la media risultante dia luogo ad un suono di una precisa tipologia (essenzialmente pitagorico oppure naturale), l'intento \`e quello di dar luogo ad un calcolo che esaurisca iterativamente tutte i possibili medi proporzionali, per conquistare alla fine una scala di suoni in buoni rapporti fra loro.

\noindent
La prima indagine che pone, comprensibilmente, ${\cal I}={\cal T}$, l'insieme dei quattro suoni pitagorici consonanti (\ref{tavoletta}), non pu\`o che confermare ci\`o che \`e emerso dalla tabella precedente, relativa a ${\mathfrak S}_{\cal P}$: posto come restrizione ${\mathfrak R}$ l'appartenenza ai suoni pitagorici (\ref{suonipit}), indicata in modo conciso dal pedice $(2,3)$, il  
calcolo (\ref{genmedie}) \`e
\begin{equation}
	\label{pitgen}
{\it\Bbb M}\left\langle {\cal T}\right\rangle_{(2,3)}=\left\{m_A(1,\frac{3}{2})=\frac{4}{3}, \;
m_A(1,2) = \frac{3}{2}\right\} \subset {\cal T}
\end{equation}
e non trova accesso ad un secondo passaggio. 
In altri termini, non esistono suoni intermedi di tipo pitagorico le cui corde formino proporzioni armoniche con le quattro prescritte. 

\noindent
Un quadro teorico pi\`u interessante, da abbinarsi alla seconda tabella attinente a ${\mathfrak S}_{\cal N}$, emerge se si allenta la restrizione ammettendo i suoni naturali (\ref{suoninat}), indicata stavolta con il pedice $(2,3,5)$: il primo passo di (\ref{genmedie}), sempre a partire dall'insieme (\ref{tavoletta}), \`e

\begin{equation}
	\label{primait}
	{\Bbb M}\left\langle {\cal T} \right\rangle_{(2,3,5)}=
	\left\{ m_A(1,\frac{3}{2})=\frac{5}{4}, \;
	m_A(1,2)=\frac{3}{2}, \; m_A(\frac{4}{3},2)=\frac{5}{3} \right\}
\end{equation}
e non sfugge il fatto che d\`a luogo all'ampliamento ${\cal T}_5$, definito in (\ref{div5}), di ${\cal T}$. Si prosegue dunque calcolando
$$
{\Bbb M}\left\langle {\cal T}_5 \right\rangle_{(2,3,5)}=
{\Bbb M}\left\langle {\cal T} \right\rangle_{(2,3,5)} \bigcup 
\left\{ m_A(1,\frac{5}{4})=\frac{9}{8}, \;
m_A(1,\frac{5}{3})=\frac{4}{3}, \; m_A(\frac{4}{3}, \frac{5}{3})=\frac{3}{2} \right\}
$$
che permette di proseguire accogliendo la novit\`a $\dfrac{9}{8}$:
$$
{\Bbb M}\left\langle {\cal T}_5 \bigcup \left\{ \dfrac{9}{8}\right\} \right\rangle_{(2,3,5)}=
{\Bbb M}\left\langle {\cal T}_5 \right\rangle_{(2,3,5)} \bigcup 
\left\{ m_A(\frac{9}{8}, 2)=\frac{25}{16}\right\}
$$
poi $\dfrac{25}{16}$:
$$
{\Bbb M}\left\langle {\cal T}_5 \bigcup \left\{ \dfrac{9}{8}, \dfrac{25}{16} \right\} \right\rangle_{(2,3,5)}=
{\Bbb M}\left\langle {\cal T}_5 \bigcup \left\{ \dfrac{9}{8} \right\} \right\rangle_{(2,3,5)} \bigcup 
\left\{ m_A(\frac{5}{4},\frac{25}{16} )=\frac{45}{32}\right\},
$$
di seguito $\dfrac{45}{32}$: 
$$
{\Bbb M}\left\langle {\cal T}_5 \bigcup \left\{ \dfrac{9}{8},  \dfrac{45}{32},\dfrac{25}{16}\right\} \right\rangle_{(2,3,5)}=
{\Bbb M}\left\langle {\cal T}_5 \bigcup \left\{ \dfrac{9}{8}, \dfrac{25}{16} \right\} \right\rangle_{(2,3,5)} \bigcup 
\left\{ m_A(\frac{9}{8},\frac{45}{32} )=\frac{81}{64}\right\}
$$
infine $\dfrac{81}{64}$:
$$
{\Bbb M}\left\langle {\cal T}_5 \bigcup \left\{ \dfrac{9}{8}, \dfrac{81}{64},  \dfrac{45}{32},\dfrac{25}{16}\right\} \right\rangle_{(2,3,5)}={\Bbb M}\left\langle {\cal T}_5 \bigcup \left\{ \dfrac{9}{8},  \dfrac{45}{32},\dfrac{25}{16}\right\} \right\rangle_{(2,3,5)}.
$$
L'ultima uguaglianza mostra che il procedimento ha avuto termine, non essendo stata generata alcuna media nuova di tipo naturale. Vale la pena elencare la gamma finale dei suoni 
\begin{equation}
	\label{sf1}
	{\mathfrak S}^{(1)}_{\cal N}=
	\bigg\{1, \frac{9}{8}, \frac{5}{4}, \frac{81}{64}, \frac{4}{3},
	\frac{45}{32}, \frac{3}{2}, \frac{25}{16}, \frac{5}{3}, 2 \bigg\}
\end{equation}
che amplia l'insieme di partenza (\ref{tav}) incorporando sei suoni. Il pedice ${\cal N}$ fa riferimento sempre ai suoni di tipo (\ref{suoninat}), l'apice $(2)$ comunica l'esistenza di (almeno) una seconda procedura iterativa che utilizza (\ref{genmedie}): effettivamente un'esplorazione del tutto spontanea viene in mente ponendo l'avvio in ${\cal I}={\mathfrak S}_N$, per rispondere alla questione: la scala naturale \`e ``chiusa'' rispetto al calcolo della media o necessita di altri suoni naturali per esaurire le medie $m_A$? Inoltre: il procedimento ha un termine come per (\ref{sf1}), cos\`\i$\,$ da dare origine ad una nuova scala di suoni naturali?
Calcoliamo dunque (radunando praticamente le caselle a sfondo grigio chiaro e grigio scuro della tabella relativa a ${\mathfrak S}_{\cal N}$)
$$
\begin{array}{l}
{\Bbb M}\left\langle {\mathfrak S}_{\cal N} \right\rangle_{(2,3,5)}=
	\left\{
	m_A(1,\frac{5}{4})=\frac{9}{8},\;
	m_A(1,\frac{3}{2})=\frac{5}{4},\;
	m_A(1,\frac{5}{3})=\frac{4}{3},\;
	m_A(1,2) = m_A(\frac{9}{8},\frac{15}{8})=m_A(\frac{4}{3},\frac{5}{3})=\frac{3}{2},\;
	\right. \\
	\\
	\left.
	m_A(\frac{9}{8},2)=
	m_A(\frac{5}{4},\frac{15}{8})=\frac{25}{16},\;
	m_A(\frac{4}{3},2)=\frac{5}{3},\;
	m_A(\frac{3}{2},\frac{15}{8})=\frac{27}{16}\;
	\right\} 
\end{array}
$$
per trovare i due suoni naturali $\frac{25}{16}$ e $\frac{27}{16}$. 
Si pu\`o procedere oltre, per effettuare tre iterazioni:
$$
{\Bbb M}\left\langle {\mathfrak S}_{\cal N} \bigcup \left\{ \dfrac{25}{16}, \dfrac{27}{16}\right\} \right\rangle_{(2,3,5)}=
{\Bbb M}\left\langle {\mathfrak S}_{\cal N} \right\rangle_{(2,3,5)} \bigcup \left\{ m_A(\frac{9}{8}, \frac{27}{16})=m_A(\frac{5}{4}, \frac{25}{16})=\dfrac{45}{32}\right\}
$$	

$$
{\Bbb M}\left\langle 	{\mathfrak S}_{\cal N}\bigcup \left\{\dfrac{45}{32}, \dfrac{25}{16}, \dfrac{27}{16}\right\} \right\rangle_{(2,3,5)}=
{\Bbb M}\left\langle 	{\mathfrak S}_{\cal N}\bigcup \left\{\dfrac{25}{16}, \dfrac{27}{16}\right\} \right\rangle_{(2,3,5)}\bigcup \left\{ m_A(\frac{9}{8}, \frac{45}{32})=\frac{81}{64}\right\}
$$	

$$
{\Bbb M}\left\langle 	{\mathfrak S}_{\cal N}\bigcup \left\{\dfrac{81}{64}, \dfrac{45}{32}, \dfrac{25}{16}, \dfrac{27}{16}\right\} \right\rangle_{(2,3,5)}=
{\Bbb M}\left\langle 	{\mathfrak S}_{\cal N}\bigcup \left\{\dfrac{45}{32}, \dfrac{25}{16}, \dfrac{27}{16}\right\} \right\rangle_{(2,3,5)}
$$
l'ultima delle quali segnala che il processo si \`e arrestato. Conviene elencare la lista finale dei suoni originali e di quelli via via aggiunti:
\begin{equation}
	\label{sf2}
	{\mathfrak S}^{(2)}_{\cal N}=
	\bigg\{1, \frac{9}{8}, \frac{5}{4}, \frac{81}{64}, \frac{4}{3},
	\frac{45}{32}, \frac{3}{2}, \frac{25}{16}, \frac{5}{3}, \frac{27}{16}, \frac{15}{8}, 2 \bigg\}.
\end{equation}
Restando nell'ordine di idee della procedura appena conclusa, osserviamo che, applicando il metodo ai suoni pitagorici ${\mathfrak S}_{\cal P}$, non viene prodotta alcuna novit\`a, in linea con (\ref{pitgen}): ovvero, ponendo ${\cal I}={\mathfrak S}_{\cal N}$ e calcolando (\ref{genmedie}) nella restrizione ${\mathfrak R}$ dei suoni pitagorici si trova 
\begin{equation}
\label{mediesp}
{\Bbb M}({\mathfrak S}_{\cal P})_{(2,3)}\subset {\mathfrak S}_{\cal P},
\end{equation}
 come mostra la tabella delle medie relativa ai suoni pitagorici ${\mathfrak S}_{\cal P}$.

\noindent
Prima di dare un'interpretazione sulla disposizione dei suoni (\ref{sf1}), (\ref{sf2}) lungo il diapason e di valutarne l'interesse musicale, prendiamo in esame la scala dei suoni equabili (\ref{scalatemp}). Le propriet\`a dapprima elencate sui medi proporzionali e la (\ref{fil}) stabiliscono la seguente relazione fra le lunghezze $a$ e $b$ delle corde e le altezze del suono $\nu_a$, $\nu_b$:
$$
m_G(\nu_a, \nu_b)=\dfrac{\kappa}{m_G(a,b)},
$$
dunque la presenza di proporzioni geometriche sulla scala dei suoni corrisponde esattamente alla presenza della proporzione geometrica delle corrispondenti lunghezze. Se si considera poi che (\ref{scalatemp}) \`e essa stessa una successione geometrica, allora \`e semplice concludere che ciascuna terna di suoni alla medesima distanza (compresi tre suoni consecutivi) forma una proporzione geometrica. L'occorrenza gi\`a incontrata (accennata sulla tavoletta) \`e l'equidistanza dei diapason disposti ad altezze in progressione geometrica. L'assenza di proporzioni armoniche nella scala equabile e il pur fitto ma disadorno intervento della proporzione geometrica nella scala equabile rende la formazione di quest'ultima non in dialogo con le altre, se non per confrontare i valori irrazionali (\ref{scalatemp}) con le frazioni naturali prossime ad essi.


\section{Sul contenuto musicale del risultato}

\noindent
Dapprima hanno recitato i personaggi (matematici), sveliamo gli interpreti (musicali) -- in parte gi\`a noti -- e la loro parte.
Innanzitutto, la capacit\`a di generare suoni che amplino la presenza della proporzione armonica nella consonanza \`e assente nelle sonorit\`a pitagoriche, come indica la (\ref{pitgen}):
la formazione di proporzioni armoniche \`e impossibile se la gamma sonora viene limitata ai suoni generati da bipartizione o tripartizione della corda.
Allo stesso tempo, la (\ref{mediesp}) testimonia che l'aggiunta dei suoni pitagorici di (\ref{scalapit}) alle quattro consonanze (\ref{tav}) non predispone alcuna proporzione armonica che intensifichi le due iniziali del diapente e del diatessaron.
In termini del processo di generazione per diapente illustrato nel Paragrafo 2.1, l'attenzione del suono aggiunto esclusivamente alla consonanza ``perfetta''\footnote{Questo \`e proprio il termine della teoria musicale per denotare le consonanze di diapente e diatessaron.} che forma con quello che lo genera 
produce un assetto finale refrattario alla disposizione armonica generale dei suoni. 
Non \`e impropria o forzata la considerazione che per vari secoli la musica monodica del canto gregoriano e l'arcaica sovrapposizione di suoni solo per diapason, diapente e diatessaron del primo Medioevo riflettono l'assetto ``ostile'' alla reciproca armonia dell'utilizzata gamma di suoni pitagorici.
L'uso dei suoni naturali (\ref{scalanat}), che adotta le consonanze ``imperfette'' $5/4$ e $5/3$, proviene appunto dall'esigenza di armonia fra suoni in un periodo in cui la polifonia e la musica strumentale prendono il sopravvento, a partire dall'alto Medioevo: se \`e la proporzione armonica a valorizzare, esaltare le consonanze, la diffusa presenza di questa nei suoni naturali ${\mathfrak S}_{\cal N}$ \`e un indicazione qualitativa della fondatezza del criterio.  

\noindent
Nondimeno, si \`e pensato di andare oltre al riscontro oggettivo offerto dai calcoli eseguiti nella tabella del Paragrafo 3.1, per affiancare al meglio la matematica delle proporzioni all'aggregazione di suoni naturali, o addirittura assegnarle un ruolo esclusivo.
Perci\`o, affidandosi al generatore di medie che seleziona solo suoni naturali ed in proporzione armonica con gli esistenti, la base delle consonanze perfette (\ref{tav})
si estende, dopo un numero finito di passi, ai suoni naturali (\ref{sf1}).
A questo punto, ci\`o che si raccoglie in (\ref{sf1}) abbandonandosi al calcolo, alla sentenza dei numeri, propone un'interpretazione dal punto di vista musicale?

\noindent
Si pu\`o cominciare col dire quello che non si \`e trovato, ovvero la scala naturale 
(\ref{scalanat}): \`e un aspetto tutt'altro che sfavorevole, 
in quanto tale scala\footnote{Nel nostro tempo nota come {\it scala di DO maggiore}.} non \`e attinente al periodo storico della teorizzazione
dei suoni naturali. Il metodo non d\`a concessioni ad
inseguire la scala dei sette suoni, volont\`a quest'ultima che trasmette probabilmente il maggior difetto nella letteratura
a proposito, che ingabbia ogni risultato nella griglia della tonalit\`a e che 
articola i passaggi nei termini delle denominazioni moderne delle note, a cui siamo abituati.
Propriamente, la raccolta (\ref{sf1}) o la (\ref{sf2}) non va intesa come sequenza di suoni da scandire in successione, una scala ascendente, ma come repertorio di moduli costruttivi in armonia tra loro, predisposti alla combinazione per la compiacenza della proporzione.
La prima evidenza apparentemente scomoda \`e l'esistenza di due suoni molto vicini, $5/4$ (il MI naturale) e $81/64$ (il MI pitagorico). In realt\`a, se inseriti nel giusto contesto, le due presenze rintracciano i due maggiori capisaldi della lunga fase di teoria e prassi musicale che precorre il Rinascimento: da una parte, i quattro suoni pitagorici

\begin{equation}
\label{finales}
\begin{array}{l}
\dfrac{9}{8}
			\hspace{1.2truecm}
			\dfrac{81}{64} \hspace{1.truecm} \dfrac{4}{3}\hspace{1.2truecm} \dfrac{3}{2}
\\
\underbrace{\hspace{1.5truecm}}_{epogdoon} \;
\underbrace{\hspace{1.3truecm}}_{limma} \;
\underbrace{\hspace{1.5truecm}}_{epogdoon}
\end{array}
\end{equation}
formano un tetracordo diatonico\footnote{Gli altri due generi, cromatico ed enarmonico furono abbandonati nel passaggio dal mondo classico a quello medievale.}, particella fondamentale della teoria armonica dell'antica Grecia, che confluisce nella teoria musicale dei modi\footnote{``Il significato del termine {\it modalit\`a} \`e talmente esteso che si qualifica come modale, con maggiore o minore propriet\`a, buona parte di ci\`o che non rientra nell'ambito del nostro sistema tonale'', \cite{utet}.} come elemento centrale per la metodica classificazione del repertorio, mediante le quattro {\it finales}, note sulle quali la melodia si riposa definitivamente e attorno alle quali si sviluppa il canto, nominate  {\it protus}, (RE, $9/8$), {\it  deuterus} (MI, $81/64$), {\it tritus} (FA, $4/3$) e {\it tetrardus} (SOL, $3/2$) e qui legittimamente presenti nella loro isolata ed antica intonazione pitagorica.
I quattro modi indotti dalle {\it finales} rappresentano il punto stabile di riferimento non solo per la sistemazione universale del canto gregoriano operata nell'VIII secolo (l'{\it octoechos} bizantino), ma persino per l'inoltrato repertorio rinascimentale di madrigali, mottetti, toccate, ... tuttora di complicato inquadramento sistematico (\cite{harold}, \cite{azzaroni}).
La compresenza in (\ref{sf1}) dell'altro MI, ad altezza $5/4$, ben documenta la dispersione dei suoni pitagorici a vantaggio di altre soluzioni e ben si inserisce nella reperibile selezione dell'{\it esacordo naturale}

\begin{equation}
\label{esnat}
\begin{array}{l}
1 \hspace{1.2truecm} \dfrac{9}{8}	\hspace{1.1truecm}
	\dfrac{5}{4} \hspace{.9truecm} \dfrac{4}{3}\hspace{1.2truecm} \dfrac{3}{2}    \hspace{1.1truecm} \dfrac{5}{3}
\\
			\\
\underbrace{\hspace{1.4truecm}}_{t.~mag.~} \;
\underbrace{\hspace{1.3truecm}}_{t.~min.~} \,
\underbrace{\hspace{.9truecm}}_{s.~mag.~}  \,
\underbrace{\hspace{1.3truecm}}_{t.~mag.} \;
\underbrace{\hspace{1.3truecm}}_{t.~min.} 
\end{array}
\end{equation}
che dalla esplicita formulazione di Guido Monaco rappresent\`o per vari secoli il sistema 
di riferimento basilare per la teoria, la pratica e la didattica musicale\footnote{\`E il cosiddetto metodo di {\it solmisazione}.}. Pur ambientandosi nella griglia della modalit\`a, l'esacordo imprime un segno fondamentale nel dirigere i tanti fermenti di un'epoca in continua evoluzione, quella tardomedievale e rinascimentale, densa di materiali sonori eterogenei e si distingue come sistema proiettato nel futuro, per esempio nella prassi di passaggio da un esacordo ad un altro pi\`u grave o pi\`u acuto (la cosiddetta {\it mutatio}), qualcosa che al giorno d'oggi \`e assimilabile alla modulazione da una tonalit\`a ad un'altra.
L'intreccio tra le finales modali e l'esacordo -- cardini per la fluida dottrina musicale di tanti secoli -- costituisce l'impianto teorico pi\`u convincente per passare in rassegna almeno un millennio di poliedrica teoria e prassi musicale.

\noindent
In (\ref{sf1}) altre presenze interessanti oltre alle quattro {\it finales} e all'esacordo 
possono essere rintracciate, per un contributo stavolta in termini di consonanze:
nell'enucleare il prospetto

$$
\begin{array}{l}
1 \hspace{1.2truecm} \dfrac{5}{4}	\hspace{1.1truecm}
\dfrac{4}{3} \hspace{.9truecm} \dfrac{3}{2}\hspace{1.2truecm} \dfrac{5}{3}    \hspace{1.1truecm} 2
\\
\\
\hspace{1.5truecm} \underbrace{\hspace{3.8truecm}\underbrace{\hspace{1.6truecm}}_{6/5}}_{8/5}
\end{array}
$$
si trova l'impronta del principio proclamato da Zarlino in \cite{zarldim} delle consonanze musicali basato sul {\it Senario}, ovvero la successione dei primi sei numeri interi che estendono il concetto pitagorico di consonanza basato sui primi quattro numeri interi. Per Zarlino la consonanza proviene dai rapporti ``superparticolari''  o ``superpartientes'' (frazioni dove il numeratore supera di un'unit\`a il denominatore)  $2/1$, $3/2$, $4/3$, $5/4$, $6/5$ e dai  
rapporti $5/3$ e $8/5$, questi ultimi collegati faticosamente al Senario, ma di cruciale importanza per l'affermarsi della musica polifonica ed il nascere della musica strumentale del tempo\footnote{Nella terminologia attuale $6/5$, $5/3$ e $8/5$ danno luogo agli intervalli di terza minore, sesta maggiore e sesta minore, rispettivamente; una linea melodica che si sovrappone ad un'altra utilizza largamente gli intervalli di terza e sesta per contrapporsi armonicamente a quella esistente.}.
Nonostante la selezione operata sulle distanze possa apparire mirata ad esaltare la conclusione, \`e semplice vedere che tutte le altre distanze -- ad eccezione di tono maggiore, tono minore e semitono maggiore, da non considerare consonanze --
riproducono sempre i medesimi valori.

\noindent
I restanti suoni $45/32$ e $25/16$ di (\ref{sf1}) mettono a disposizione una  proposta in merito all'articolata questione della divisione del tono: a partire dal limma pitagorico\footnote{Elencare le altre sottounit\`a di genere pitagorico, calcolabili sulla base della prosecuzione del ciclo dei diapente, sposterebbe l'attenzione su altri argomenti.} contrassegnato dall'impervio rapporto $256/243$, il problema di infittire ulteriormente il diapason ha dato vita nei secoli ad un'appassionata ricerca che ha stimolato i teorici anche dal punto di vista del calcolo aritmetico e della geometria applicati all'acustica\footnote{Fra i numerosi casi, nominiamo il trattato {\it Alia Musica}, di provenienza francese o fiamminga, in cui compare la soluzione geometrica per il calcolo del medio proporzionale armonico.}: ad esempio
\`e di Zarlino stesso la dimostrazione che il tono maggiore non pu\`o essere suddiviso in due intervalli uguali, ovvero, in termini matematici, la frazione $9/8$ non \`e il prodotto di  numeri razionali uguali.
In merito alla ricerca di gradini pi\`u piccoli che siano in armonia con il resto, la risoluzione (\ref{sf1}) addita indicazioni nei suoni dell'intervallo fra $4/3$ e $5/3$, qui rappresentati\footnote{Se pu\`o orientare un'odierna denominazione di massima, esulando la questione dei semitoni diatonici o cromatici, $45/32$ compreso tra $4/3$ (FA) e $3/2$ (SOL) corrisponde al FA diesis, $25/16$ tra $3/2$ (SOL) e $5/3$ (LA) corrisponde al SOL diesis.}

$$
\begin{array}{l}
\overbrace{\hspace{3.1truecm}}^{tono\;maggiore} \;\;
\overbrace{\hspace{2.9truecm}}^{tono\;minore} \\
\dfrac{4}{3} \hspace{1.2truecm} \dfrac{45}{32}	\hspace{1.3truecm}
\dfrac{3}{2} \hspace{1.truecm} \dfrac{25}{16}\hspace{1.3truecm} \dfrac{5}{3}   
\\
\\
\underbrace{\hspace{1.5truecm}}_{135/128} \;
\underbrace{\hspace{1.6truecm}}_{16/15} \;
\underbrace{\hspace{1.3truecm}}_{25/24} \;
\underbrace{\hspace{1.5truecm}}_{16/15} 
\end{array}
$$
Effettivamente nella suddivisione sia del tono maggiore che del tono minore \`e evidente il ruolo del semitono maggiore di misura $16/15$, gi\`a conosciuto nel disporre i suoni naturali (\ref{scalanat}). D'altra parte, anche il frazionamento  $25/24$, il cosiddetto {\it semitono minore}, altro elemento centrale nelle suddivisioni armoniche di Zarlino, \`e importante per gestire ed arricchire il repertorio dei suoni.
L'inevitabile sbalzo $135/128$ che occorre per ricoprire il tono maggiore non \`e 
privo di significato: cogliamo l'occasione per scrivere $\frac{135}{128}=\frac{5}{4}\times (\frac{3}{2})^2 \times \frac{1}{2^3}$ e tradurre l'operazione come ``aggiunta di due diapente (o quinte) al suono $5/4$ (MI) e trasporto in grave di due diapason (o ottave)'', con la mera intenzione di un unico accenno 
al modo di ragionare (in tutti i tempi) sul materiale sonoro mediante le frazioni, abbracciando tuttavia contenuti e metodi che vanno ben oltre i nostri intenti\footnote{Per chi ha la pazienza di seguire la soluzione, semplificando in suoni equabili: tre quinte dal MI fa percorrere MI$\,\rightarrow \, SI\,\rightarrow\,FA\#\,\rightarrow\,DO\#$, quest'ultimo da trasportare due ottave sotto, per trovare un suono fra il DO ed il RE della scala originale.}.

\noindent
Rivolgiamo ora l'attenzione alla seconda lista (\ref{sf2}), frutto dell'esperimento di ``chiusura'' della scala naturale (\ref{scalanat}) tramite il generatore di medie: ovviamente $\mathfrak{S}_{\cal N}^{(2)}$ amplia $\mathfrak{S}_{\cal N}^{(1)}$, coerentemente alla posizione ${\cal T}\subset \mathfrak{S}_{\cal N}$ degli insiemi di partenza, ma non di tanto: le uniche nuove presenze sono $\frac{27}{16}$ e $\frac{15}{8}$. 
I due suoni\footnote{Essi corrispondono rispettivamente al LA pitagorico e al SI di intonazione naturale.} vanno a riempire lo spazio vuoto pi\`u ampio di ${\mathfrak S}_{\cal N}^{(1)}$, quello tra $5/3$ e $2$:
$$
\begin{array}{l}
\overbrace{\hspace{3.truecm}}^{tono\;maggiore} \\
\hspace{1.truecm}\overbrace{\hspace{2.truecm}}^{tono\;minore} \\
\dfrac{5}{3} \hspace{0.4truecm} \dfrac{27}{16}	\hspace{1.5truecm}
\dfrac{15}{8} \hspace{1.truecm} 2   
\\
\underbrace{\hspace{1.truecm}}_{81/80} \hspace{1.2truecm}
\underbrace{\hspace{1.5truecm}}_{semitono\;maggiore} 
\end{array}
$$
senza tuttavia introdurre novit\`a dal punto di vista della misura degli intervalli.
La particella $81/80$ \`e il pezzo pi\`u piccolo della collezione: fa parte dei cosiddetti {\it comma} -- termine traducibile appunto con ``frammento'' -- ed \`e talmente ricorrente nella trattatistica musicale da ricevere svariate denominazioni\footnote{Le principali sono comma diatonico, sintonico, didimeo, tolemaico.} a seconda della mansione svolta: qui \`e evidente il ruolo di scarto tra il tono maggiore e tono minore, come d'altra parte lo si pu\`o ottenere dal divario dei due MI pitagorico e naturale dapprima considerati: $\frac{81}{64}: \frac{5}{4}=\frac{81}{80}$. 
La comparsa del suono pitagorico $\frac{27}{16}$ si collega ad almeno due chiavi di lettura interessanti: da una parte, si struttura un nuovo impianto quattro finales--esacordo nei valori, rispettivamente, $\left( \frac{27}{16}, \frac{15}{8}, 2, \frac{9}{4}\right)$ e $\left( \frac{3}{2}, \frac{5}{3}, \frac{15}{8}, 2, \frac{9}{4}, \frac{5}{2} \right)$ -- considerando al di l\`a del $2$ i suoni omologhi sul diapason pi\`u acuto --
con il medesimo dualismo dei due suoni addossati, l'uno $\frac{27}{16}$ pitagorico l'altro $\frac{5}{3}$ naturale, ma con distanze reciproche disposte differentemente rispetto alla precedente presentazione (\ref{finales}), (\ref{esnat})\footnote{Si tratta stavolta del cosiddetto {\it esacordo duro}.} e con un colore diverso. Dall'altra parte, il suono $\frac{27}{16}$ va a completare la presenza di tutti i diapente dell'esacordo naturale, nel senso che tutti i suoni di (\ref{esnat}) trasportati in alto di un diapente (ovvero moltiplicati per $3/2$) tornano su un suono dell'esacordo medesimo, a meno eventualmente di un irrilevante diapason, eccetto $\frac{9}{8}$, che appunto si trasferisce su $\frac{9}{8}\times \frac{3}{2}=\frac{27}{16}$.

\noindent
Malgrado la non vistosa differenza fra i due insiemi (\ref{sf1}) e (\ref{sf2}), un fatto \`e da segnalare: a partire dal quadricordo ${\cal T}$ di consonanze essenziali non \`e possibile ottenere le sette note diatoniche formanti, piegandosi al linguaggio della tonalit\`a, la ``scala di DO maggiore''; il settimo suono SI deve essere gi\`a presente, come avviene in (\ref{sf2}), affinch\'e possa essere annoverato.
La ``scala'' apparentemente pi\`u consona alla disposizione armoniosa dei suoni -- ovviamente in base al punto di vista percorso -- sembra essere quella esacordale, che si afferma in (\ref{sf1}) insieme ai frazionamenti dapprima esaminati.

\noindent
Commentiamo infine la sistemazione equidistante (\ref{scalatemp}),
al cospetto delle proporzioni: l'unica presente \`e quella geometrica, in modo peraltro abbondante, dal momento che suoni a medesima distanza da un suono centrale formano proporzioni geometriche.
Dal punto di vista dell'eufonia -- e ci\`o ribadisce la fiducia riposta nella proporzione armonica -- l'accostamento dei tre suoni in tale proporzione \`e sostanzialmente privo di significato musicale: sono in proporzione geometrica ad esempio tre suoni intervallati da due diapason, perfettamente omogenei, come tre suoni emessi da tre tasti consecutivi di un pianoforte, di ascolto non certo gradevole.
La questione della consonanza nel temperamento equabile, designato come quello correntemente predominante, almeno in teoria, \`e molto complessa
e coinvolge argomenti pi\`u assimilabili alla simmetria che alla proporzione; l'ambito musicale pi\`u frequentato da questa analisi \`e il repertorio otto-- novecentesco, lontano dal periodo in cui ci siamo calati: un esempio in tal senso \`e l'articolo \cite{simm}.

\section{Coda}

\noindent
Le motivazioni informali e personali del lavoro ruotano attorno all'interrogativo ``come mai \`e cos\`\i$\,$ immediata, trasparente la formazione, la struttura matematica della scala pitagorica e della scala equabile -- connaturate rispettivamente  con la generazione per diapente in successione ed con l'equidistanza attraverso il diapason --  quanto invece \`e apparentemente assente, introvabile nella scala naturale, presentata sempre -- nei testi di teoria musicale -- in modo manuale, enumerando l'elenco dei suoni?''

\noindent
Probabilmente la bellezza dei suoni naturali sta proprio nella difficolt\`a di 
far emergere il ruolo esclusivo della matematica che chiarisca, nella vasta scelta (\ref{suoninat}), la sequenza di Zarlino: la matematica armonia del calcolo avviato in (\ref{riemp}) compie un tentativo in tal senso.

\noindent
D'altra parte, pu\`o essere la matematica stessa in veste armonica a cercare con i propri strumenti i suoni: 
il modesto conteggio sviluppato nel lavoro ha voluto affiancare l'idea che un incastro di proporzioni armoniche possa assegnare ai suoni naturali un'anima razionale, teorizzabile. 
A partire dall'ossatura inconfutabile della cetra dei quattro suoni pitagorici (\ref{tav}), un percorso in tal senso non pu\`o non incontrare il principio maestro nel Rinascimento di proporzione metrica, per guidare l'armonia dei componenti, e dei componenti con il tutto.
\`E per questo comprensibile, immagino, l'estemporaneo abbandono in qualche punto di ci\`o che precede alla sconfinata ammirazione per capolavori dell'architettura, qui trascinati dalla musica. 
Gli esiti raccolti in (\ref{sf1}) e (\ref{sf2}) propongono uno scrigno sonoro
di tutto rispetto, fedele al percorso storico--musicale fino al 500, atto se non altro per accendere la riflessione sui veri impianti teorici che sorreggono la musica (almeno) fino al tempo di Zarlino e scansare la troppo frequente forzatura atta ad incasellare i suoni nelle sette note della scala diatonica (i tasti bianchi del pianoforte) o, peggio ancora, nelle dodici note della scala cromatica (tasti bianchi e neri). L'emerso predominio dell'esacordo, con la sola cifra dell'algoritmo, \`e una rivincita sulla leggerezza di sentir dire spesso che ``Guido Monaco ha inventato le sette note''.

\noindent
Al di l\`a del significato musicale proposto ed analizzato -- peraltro in modo onesto e non forzato, anche se consapevolmente impugnabile con svariati argomenti -- c'\`e un aspetto che ritengo degno di interesse: 
la procedura \`e chiusa, ha avuto un inizio ed una fine; il progetto di costruzione si \`e esaurito, \`e approdato ad una conclusione. La contrapposizione spontanea \`e offerta dagli infiniti suoni pitagorici generati per diapente, ciascuno dei quali rende conto della propria consonanza solo al precedente ma non tiene memoria di tutti gli altri, con i quali non \`e tenuto a formare proporzioni sonore. L'istintivo richiamo a valutazioni di carattere estetico e concettuale dei nostri tempi\footnote{ Nella vasta letteratura a proposito citiamo Cesare Brandi, 1906--1988, illustre storico dell'arte e critico d'arte: 
``{\it La differenza cruciale tra gotico e Rinascimento sta tutta qui: la cancellazione dell'infinito come slancio mistico e dell'infinito come dimensione in cui si sciolgono le guglie gotiche. Lo spazio si rinchiude dentro l'uomo.''}.}
riguardo al lungo periodo culturale preso in esame sembra avvalorare l'idea di una torre infinita di suoni medievali, dove ciascun pezzo riconosce solo quello sottostante, contrapposta al disegno compiuto di suoni rinascimentali, dove gli elementi si riconoscono e si completano nelle proporzioni.
Tentando di esprimere una percezione sfuggente, non semplice da mettere a fuoco, direi che il trovato compimento \`e proprio espressione del rispetto imposto dalla proporzione, dell'armonia della parte con il tutto, in accordo con l'idea rinascimentale, albertiana di perfezione nel senso etimologico di {\it perfectum}, ultimato, concluso.

\noindent
Proseguendo il capriccioso parallelo tra i mondi sonori esplorati e la pluralit\`a di invenzioni dell'architettura -- per attardarsi nel tema dominante dell'articolo -- la sublime levatura dei ristretti numeri (\ref{tavoletta}) \`e paragonabile alla perfezione di poche linee mirabilmente combinate di un tempio greco, in cui il numero \`e al servizio del sacro; la fredda e controllata disposizione dei suoni equabili fa pensare alla regolare ripetitivit\`a di un grattacielo moderno, in cui il numero \`e assoggettato al pratico, al funzionale, sempre uguale a se stesso in una parte o nella vista di insieme.
Non \`e cos\`\i$\,$ per la percezione di unit\`a in cui si raccolgono le parti in rapporto armonico, per l'inesauribile ricerca di nuovi spunti e dettagli che si crea attorno al gi\`a formato che ci permea quando si \`e immersi nei capolavori architettonici del Rinascimento animati dalle proporzioni, le medesime che troviamo nei suoni naturali.

\noindent
\`E un rapporto difficile quello tra matematica e musica, \`e piena di insidie l'esplorazione oggettiva senza essere manovrati dalla suggestione di rintracciare la matematica ovunque, come avviene troppo spesso nell'analisi formale e strutturale del repertorio musicale.
Sono poche le leve scientifiche che lanciano un'ispirazione
emotiva, una percezione estetica verso le categorie dell'assoluto, per mezzo del riscontro della razionalit\`a: tra le poche, il principio architettonico--musicale del Rinascimento ha probabilmente il posto d'onore.

\end{document}